\documentclass[10pt]{amsart}
\usepackage[margin=3cm]{geometry}
\usepackage{graphicx} 
\usepackage{subcaption}

\usepackage[T1]{fontenc}
\usepackage{amssymb}            
\usepackage{mathtools}    
\usepackage{mathrsfs}
\DeclareSymbolFontAlphabet{\mathrsfs}{rsfs}
\usepackage[mathscr]{eucal}    
\usepackage{tikz}              
  \usetikzlibrary{cd}          

\usepackage{booktabs}
\usepackage{array}

\usepackage[style=alphabetic]{biblatex}
\usepackage{hyperref}
\hypersetup{
    colorlinks=true,
    linkcolor=blue,
    filecolor=magenta,
    urlcolor=cyan,
    pdftitle={Higher order discrete Poincare and Bogovskii operators},
    }
\numberwithin{equation}{section}

\usepackage{aliascnt}

\theoremstyle{plain}
\newtheorem{theorem}{Theorem}[section]

\newaliascnt{corollary}{theorem}
\newtheorem{corollary}[corollary]{Corollary}
\aliascntresetthe{corollary}

\newaliascnt{lemma}{theorem}
\newtheorem{lemma}[lemma]{Lemma}
\aliascntresetthe{lemma}

\newaliascnt{proposition}{theorem}
\newtheorem{proposition}[proposition]{Proposition}
\aliascntresetthe{proposition}

\newaliascnt{observation}{theorem}

\aliascntresetthe{observation}

\newaliascnt{conjecture}{theorem}

\aliascntresetthe{conjecture}

\theoremstyle{definition}

\newaliascnt{definition}{theorem}
\newtheorem{definition}[definition]{Definition}
\aliascntresetthe{definition}

\newaliascnt{example}{theorem}
\newtheorem{example}[example]{Example}
\aliascntresetthe{example}

\newaliascnt{notation}{theorem}

\aliascntresetthe{notation}

\theoremstyle{remark}

\newaliascnt{remark}{theorem}
\newtheorem{remark}[remark]{Remark}
\aliascntresetthe{remark}

\usepackage[nameinlink, capitalize, noabbrev]{cleveref}

\crefname{theorem}{Theorem}{Theorems}
\Crefname{theorem}{Theorem}{Theorems}

\crefname{corollary}{Corollary}{Corollaries}
\Crefname{corollary}{Corollary}{Corollaries}

\crefname{lemma}{Lemma}{Lemmas}
\Crefname{lemma}{Lemma}{Lemmas}

\crefname{proposition}{Proposition}{Propositions}
\Crefname{proposition}{Proposition}{Propositions}

\crefname{observation}{Observation}{Observations}
\Crefname{observation}{Observation}{Observations}

\crefname{conjecture}{Conjecture}{Conjectures}
\Crefname{conjecture}{Conjecture}{Conjectures}

\crefname{definition}{Definition}{Definitions}
\Crefname{definition}{Definition}{Definitions}

\crefname{example}{Example}{Examples}
\Crefname{example}{Example}{Examples}

\crefname{notation}{Notation}{Notations}
\Crefname{notation}{Notation}{Notations}

\crefname{remark}{Remark}{Remarks}
\Crefname{remark}{Remark}{Remarks}

\newcommand{\RR}{\ensuremath{\mathbb{R}}}  
\newcommand{\ZZ}{\ensuremath{\mathbb{Z}}}  

\DeclareMathOperator{\id}{id}
\DeclareMathOperator{\im}{im}
\DeclareMathOperator{\tr}{tr}       
\DeclareMathOperator{\grad}{grad}
\DeclareMathOperator{\curl}{curl}
\let\div\relax
\DeclareMathOperator{\div}{div}       
                    
\newcommand{\sspan}{\operatorname{span}}
\newcommand{\cone}{{\mathbf{P}}}         
\newcommand{\contract}{{\lrcorner\,}}      
\newcommand{\ext}{{\mathrm{ext}}}           

\newcommand{\pr}{{\mathrm{pr}}}     

\DeclarePairedDelimiter{\set}{\lbrace}{\rbrace}
\DeclarePairedDelimiter{\ipdelim}{\langle}{\rangle}
\DeclarePairedDelimiter{\abs}{\lvert}{\rvert}
\DeclarePairedDelimiter{\normdelim}{\lVert}{\rVert}

\allowdisplaybreaks

\begin{document}

\title[Higher order discrete Poincar\'e and Bogovski\u{\i} operators]{Discrete Poincar\'e and Bogovski\u{\i} operators for higher order finite element de Rham complexes}

\author{Johnny Guzm\'an}
\address{Division of Applied Mathematics,
Brown University,
Box F,
182 George Street,
Providence, RI 02912, USA}
\email{\href{mailto:johnny_guzman@brown.edu}{johnny\_guzman@brown.edu}}

\author{Pratyush Potu}
\address{Division of Applied Mathematics,
Brown University,
Box F,
182 George Street,
Providence, RI 02912, USA}
\email{\href{mailto:pratyush_potu@brown.edu}{pratyush\_potu@brown.edu}}

\author{K. Olav Tyssvang}
\address{Mathematical Institute, 
University of Oxford,
UK}
\email{\href{mailto:olav.tyssvang@linacre.ox.ac.uk}{olav.tyssvang@linacre.ox.ac.uk}}

\begin{abstract}
We present discrete Poincar\'e and Bogovski\u{\i} operators for the common finite element discretizations of the de Rham complex of arbitrary polynomial order. We first discuss how such operators arise from a decomposition of the finite element spaces into a lowest order component and a remainder. This results in a fast algorithm for computing $\grad$, $\curl$, and $\div$ potentials. We then introduce a nonlocal Poincar\'e operator on contractible domains, which generalizes existing constructions based on cone operators to arbitrary order by integrating over cones against extended test forms. This requires developing commuting extension operators, which turn out to be related to rational forms. We also show how the nonlocal construction is modified to obtain a discrete Bogovski\u{\i} operator preserving homogeneous boundary conditions. We further discuss how the ideas of the construction fit into the theory of finite element systems.
\end{abstract}

\keywords{finite element exterior calculus, Poincar\'e operator, Bogovski\u{\i}
operator, polynomial differential forms, finite element systems}
\subjclass[2020]{Primary 65N30; Secondary 58A12}

\maketitle



\section{Introduction}
The classical Poincar\'e lemma states that every closed differential $k$-form, $k\geq1$, on a star-shaped domain is exact. A standard proof uses a Poincar\'e operator. If $(\Lambda^\bullet(\Omega),d^\bullet)$ denotes the complex of smooth differential forms on such a domain $\Omega$, a Poincar\'e operator is a graded map $P^\bullet$ with $P^k:\Lambda^k(\Omega)\to\Lambda^{k-1}(\Omega)$ and
\begin{equation}\label{eq:intro--homotopy}
d^{k-1} P^k+P^{k+1}d^k =\id-\pi^k,
\qquad 0\leq k\leq n,
\end{equation}
where $P^0=P^{n+1}=0$, $\pi^k=0$ for $k>0$, and $\pi^0$ is a projection onto the constant $0$-forms. The Poincar\'e operator is a special case of a homotopy operator in differential geometry (see e.g. \cite[Section 17]{Lee:13}) and \eqref{eq:intro--homotopy} is referred to as the homotopy identity. A Bogovski\u{\i} operator is the analogue of a Poincar\'e operator for forms satisfying homogeneous boundary conditions following \cite{Bogovskii:79}.

Poincar\'e and Bogovski\u{\i} operators that map polynomial-valued differential forms to polynomial-valued differential forms have become fundamental tools in the construction and analysis of finite elements, and in particular in the finite element exterior calculus (FEEC) \cite{Arnold-Falk-Winther:06}. In \cite{Hiptmair:99}, for instance, a Poincar\'e operator restricted to simplices, under the name ``potential mapping'', is used to define finite element spaces equivalent to the trimmed polynomial spaces of \cite{Arnold-Falk-Winther:06}. Furthermore, in \cite{Costabel-McIntosh:10}, regularized Poincar\'e operators which still preserve polynomials are used to establish Poincar\'e inequalities. In general, such operators recur throughout the analysis of finite element methods \cite{Boffi-Costabel-Dauge-Demkowicz-Hiptmair:11, Braess-Pillwein-Schoberl:09, Farrell-Mitchell-Scott-Wechsung:22, Melenk-Rojik:20, arnold-guzman:21, Cap-Hu:23}.

The well-known de Rham theorem (see e.g.\ \cite[Theorem~18.14]{Lee:13}) asserts that the de Rham cohomology of a smooth manifold is isomorphic to the simplicial cohomology of any triangulation of the manifold. So, a natural question that arises is whether one can define a \emph{discrete Poincar\'e operator} directly over a simplicial cochain complex, i.e. an operator which maps $k$-cochains to $(k-1)$-cochains and satisfies the homotopy identity, without invoking the continuous level operator. The existence of such an operator is used widely in algebraic topology (see e.g.\ \cite{Munkres:84, Hatcher:02}). Constructive (i.e. computable) constructions of such operators are discussed at length in \cite{guzman-potu:26}. In particular, these constructions revolve around dualizing a ``generalized cone operator'' which maps $k$-chains to $(k+1)$-chains and satisfies the homotopy identity on the chain complex. Due to the isomorphism between cochains and Whitney forms \cite[Chapter IV, Section~27]{Whitney:57}, having such a discrete Poincar\'e operator over cochains is equivalent to having one over the Whitney forms.

Higher order generalizations of the Whitney forms are at the core of FEEC \cite{Arnold-Falk-Winther:06}. We would therefore like to know whether the discrete Poincar\'e operators defined for Whitney forms extend to higher order spaces. In this work, we answer this question in the affirmative by providing explicit discrete Poincar\'e and Bogovski\u{\i} operators for the higher order spaces. We consider in particular the complex of trimmed polynomial forms on a simplicial mesh $\mathcal T_h$, $\bigl(\mathcal P_r^-\Lambda^k(\mathcal T_h),d^k\bigr)_{k=0}^n$, and the full degree reducing complex of polynomial forms $\bigl(\mathcal P_{r-k}\Lambda^k(\mathcal T_h),d^k\bigr)_{k=0}^n$ (see \cite{Arnold-Falk-Winther:06} or \cref{section:preliminaries} for precise definitions). We present two constructions yielding Poincar\'e operators for these complexes.

In fact, such a higher order discrete Poincar\'e operator is immediately obtained from a ``discrete $\mathfrak{B}^*$ problem'' posed globally over a higher order space (see Section \ref{section:global-construction}). In \cite[Section 5]{DiPietro-Droniou-Pitassi:25}, the authors develop an alternative discrete Poincar\'e operator by decomposing the higher order spaces into the complex of Whitney forms and a remainder. Remarkably, the resulting operator is local in the sense that, once that lowest order operator has been supplied, the remaining work is cellwise and can be implemented by matrices computed once on a reference simplex. We review the construction of \cite{DiPietro-Droniou-Pitassi:25} and show an alternative way of handling the remainder via discrete $\mathfrak B^*$ problems posed over the remainder spaces. We implement the operator and provide numerical experiments showing the efficiency of the operator compared with the global problem.

In \cite{guzman-potu:26}, a nonlocal lowest order discrete Poincar\'e is developed by considering ``cones'' of simplices. In the case of a star-shaped domain, the cone is obtained as the convex hull of the simplex and some fixed point $p_0$ in the mesh. To extend this construction to the higher order spaces, we develop formulas involving integration against extensions of polynomials on simplices to cones. The necessary properties on these extensions lead to novel connections to rational forms. We then provide a generalization of this operator to more general domains by replacing the cone with the ``extrusion''.  We also show how to modify the operator to accommodate homogeneous boundary conditions, thus yielding a discrete Bogovski\u{\i} operator.

The presented discrete Poincar\'e and Bogovski\u{\i} operators have natural applications. Firstly, they provide a numerical method for computing discrete potentials: the problem of finding $u$ such that $du = f$ for a given $f$ such that $df=0$. Moreover, in the work \cite{boon:26}, it is shown that bounded discrete Poincar\'e operators allow for efficient solvers for Hodge-Laplace problems as well as auxiliary space preconditioners. 

The paper is organized as follows. \Cref{section:preliminaries} fixes the finite element spaces and notation. In \Cref{section:global-construction}, we present how a discrete Poincar\'e operator is obtained from solving a $\mathfrak{B}^*$ problem over the global discrete spaces. Then, in \cref{section:local-construction}, we show how decomposing a higher order polynomial $k$-form into its Whitney component and an exact remainder leads to a more efficient discrete Poincar\'e operator. In \Cref{section:cone-extension}, we provide the nonlocal construction based on polynomial extensions over parameterized cones. In \cref{section:bogovskii}, we present modified versions of these two constructions yielding discrete Bogovski\u{\i} operators. Finally, \Cref{section:implementation} describes the implementation and numerical experiments of the local method, and \cref{section:fes} formulates the constructions in terms of finite element systems.


\section{Preliminaries}
\label{section:preliminaries}

Throughout, $\Omega\subset\RR^n$ is a connected open bounded polytopal domain with Lipschitz boundary, oriented by the standard orientation of $\RR^n$. We assume familiarity with differential forms as used in the finite element exterior calculus (FEEC). Standard references are \cite{Arnold-Falk-Winther:06,Arnold:18}.

To fix notation, we write $\Lambda^k(D)$ for the space of smooth differential $k$-forms on a domain $D$, and $L^2\Lambda^k(D)$ for the $k$-forms with $L^2$-regular coefficients. The exterior derivative is denoted $d$ in every degree, the degree in which it acts being clear from context, and $\wedge:\Lambda^k\times\Lambda^\ell\to\Lambda^{k+\ell}$ is the wedge product. For a sub-manifold $\tau\subset D$ we denote by $\tr_\tau\coloneqq i_\tau^*$ the trace map $\Lambda^k(D)\to \Lambda^k(\tau)$, which is given by the pullback along the inclusion $i_\tau:\tau \hookrightarrow D$. The Euclidean metric on $\Omega$ induces a pointwise inner product on forms, denoted $\ipdelim{\cdot,\cdot}$, and the corresponding $L^2$ inner product is
\begin{equation}
(u,v)_\Omega\coloneqq \int_\Omega\ipdelim{u,v}\operatorname{vol}_\Omega,
\end{equation}
where $\operatorname{vol}_\Omega$ is the induced volume form on $\Omega$. For a $k$-form $u$ we denote by $|u|\coloneqq \langle u, u\rangle^{1/2}$ the pointwise absolute value. Lastly, the Sobolev spaces of differentiable forms (see \cite{Goldshtein1982}) are given by
\begin{equation}
\mathrsfs W_{p,q}^k(\Omega) \coloneqq \set{ u\in L^p\Lambda^k(\Omega): du\in L^q\Lambda^{k+1}(\Omega) },
\end{equation}
where $d$ is acting in the weak sense. In particular we are interested in the spaces $H\Lambda^k(\Omega)\coloneqq \mathrsfs W_{2, 2}^k(\Omega)$, since for $n=3$ they correspond to the physically relevant spaces $H^1(\Omega)$, $H(\curl, \Omega)$, $H(\div, \Omega)$, and $L^2(\Omega)$ for $k=0,1,2,3$ respectively.

\subsection{The mesh and its chain complex}
\label{section:prelim--mesh}

Let $\mathcal T_h$ be a finite, pure $n$-dimensional geometric simplicial complex whose underlying space is $\overline\Omega$, the closure of $\Omega$. This is exactly a conforming simplicial mesh in the sense of \cite[Definition 1.55]{Ern-Guermond:04}. Let $\mathcal T_h^\ell\subset \mathcal T_h$ denote its $\ell$-dimensional simplices, so that $\mathcal T_h=\bigcup_{\ell=0}^n\mathcal T_h^\ell$. For $\sigma \in \mathcal T_h$, we write $\tau\triangleleft\sigma$ when $\tau$ is a face (that is, a subsimplex) of $\sigma$. The case $\tau=\sigma$ is allowed, and we call $\tau$ a proper face when $\tau\neq \sigma$. We also write $\partial\mathcal T_h^\ell\subset \mathcal T_h^\ell$ for those $\ell$-simplices contained in $\partial\Omega$ and $\partial\mathcal T_h = \bigcup_{\ell=0}^{n-1} \partial\mathcal T_h^\ell$. We fix a global ordering of the vertices of $\mathcal T_h$. If $x_0<\cdots<x_\ell$ are the vertices of $\sigma \in \mathcal T_h^\ell$, we orient $\sigma$ by its increasing vertex enumeration order and write $\sigma = [x_0,\dots,x_\ell]$ to identify $\sigma$ as an oriented geometric simplex corresponding to the closed convex hull of its vertices.

The construction of \cref{section:cone-extension} invokes simplicial chains, so we fix that language now. Let $\mathcal C_\ell(\mathcal T_h)$ be the real vector space with basis $\mathcal T_h^\ell$, whose elements are the simplicial $\ell$-chains, and let $\partial_\ell:\mathcal C_\ell(\mathcal T_h)\to\mathcal C_{\ell-1}(\mathcal T_h)$ be the boundary operator,
\begin{equation}
\label{eq:prelim--boundary}
\partial_\ell[x_0,\dots,x_\ell]=\sum_{j=0}^\ell(-1)^j\,[x_0,\dots,\widehat{x_j},\dots,x_\ell]
=\sum_{\tau\triangleleft\sigma}o_{\sigma\tau}\,\tau,
\end{equation}
which satisfies $\partial_{\ell-1} \circ \partial_\ell=0$ and defines the incidence coefficients $o_{\sigma\tau}\in\set{0, 1,-1}$. Here $o_{\sigma\tau}$ is taken to be zero if $\dim \tau \ne \dim \sigma-1$. The dual space of $\mathcal C_\ell(\mathcal T_h)$ is the space of $\ell$-cochains $\mathcal C^\ell(\mathcal T_h)$. The coboundary operator $\mathsf d^\ell:\mathcal C^\ell(\mathcal T_h)\to\mathcal C^{\ell+1}(\mathcal T_h)$ is defined by duality as
\begin{equation}
(\mathsf d^\ell \mathsf{z})(c)\coloneqq \mathsf{z}(\partial_{\ell+1}c),
\qquad \mathsf{z}\in\mathcal C^\ell(\mathcal T_h),\
c\in\mathcal C_{\ell+1}(\mathcal T_h).
\end{equation}
By the de Rham theorem \cite[Theorem 18.14]{Lee:13}, the cohomology of $(\mathcal C^\bullet(\mathcal T_h),\mathsf d^\bullet)$ is isomorphic to the de Rham cohomology $H^\bullet(\Omega)$.

Integration of forms on (oriented) simplices is defined in the natural way via the pullback of the affine mapping between a fixed reference simplex and any physical simplex. This notion extends to integration of forms on chains by linearity: for $c=\sum_i a_i\sigma_i\in\mathcal C_\ell(\mathcal T_h)$ and $u$ a sufficiently regular $\ell$-form,
\begin{equation}
\label{eq:prelim--chain-integration}
\int_c u \coloneqq \sum_i a_i\int_{\sigma_i}\tr_{\sigma_i} u.
\end{equation}

\subsection{Spaces of polynomial differential forms}
\label{section:prelim--polynomial-spaces}

Finite element methods assemble a global space on $\mathcal T_h$ by gluing local spaces defined on the individual cells. The essential FEEC spaces are the two families of polynomial forms $\mathcal P_r\Lambda^k$ and $\mathcal P_r^-\Lambda^k$, which we describe now. They correspond to the classic Lagrange, Raviart-Thomas, Brezzi-Douglas-Marini, and N\'ed\'elec finite elements \cite{Arnold-Falk-Winther:06}.

Each $\sigma=[x_0,\dots,x_\ell]\in\mathcal T_h$ carries affine barycentric coordinates $(\lambda_i)_{0\le i\le \ell}$, which are scalar functions on $\sigma$ such that $\lambda_i(x_j)=\delta_{ij}$ and $\sum_i\lambda_i\equiv 1$. Fix an $\ell$-simplex $\sigma$ and an integer $r\geq0$. At each point $x\in\sigma$ the covectors $d\lambda_{i_1}\wedge\cdots\wedge d\lambda_{i_k}$, for $1\le i_1<\cdots<i_k\le \ell$, form a basis of the space $\Lambda^k\mathrm T^*_x\sigma$ of alternating $k$-linear forms on the tangent space $T_x\sigma$. For a multi-index $\alpha\in\ZZ_{\ge0}^{\ell+1}$, write $\lambda^\alpha\coloneqq\lambda_0^{\alpha_0}\cdots\lambda_\ell^{\alpha_\ell}$ and $\abs{\alpha}\coloneqq\sum_{i=0}^\ell\alpha_i$. The space of polynomial $k$-forms of degree at most $r$ is
\begin{equation}
\label{eq:prelim--full-polynomial}
\mathcal P_r\Lambda^k(\sigma)\coloneqq
\sspan\set[\big]{\lambda^\alpha\, d\lambda_{i_1}\wedge\cdots\wedge d\lambda_{i_k}
\ :\ \abs{\alpha}\le r,\ 1\le i_1<\cdots<i_k\le \ell}.
\end{equation}

A second family of polynomial forms is obtained by trimming $\mathcal P_r\Lambda^k(\sigma)$. Fix a point $a$ in $\sigma$ and let $\kappa_a \coloneqq (x-a)\contract\ :\ \Lambda^k(\sigma)\to\Lambda^{k-1}(\sigma)$ denote the associated Koszul operator. That is, the contraction against the radial vector field $X(x)= x-a$. The trimmed polynomial forms are
\begin{equation}
\label{eq:prelim--trimmed-polynomial}
\mathcal P_r^-\Lambda^k(\sigma)\coloneqq \mathcal P_{r-1}\Lambda^k(\sigma)+\kappa_a\,\mathcal P_{r-1}\Lambda^{k+1}(\sigma).
\end{equation}
This is a subspace of $\mathcal P_r\Lambda^k(\sigma)$, and does not depend on the base point $a$. The sum in \eqref{eq:prelim--trimmed-polynomial} becomes direct if $\mathcal P_{r-1}\Lambda^{k+1}$ is replaced by the subspace of forms whose coefficients are homogeneous polynomials in $\lambda_1,\cdots, \lambda_\ell$ of degree exactly $r-1$.

In degree zero one has $\mathcal P^-_r\Lambda^0(\sigma)=\mathcal P_r\Lambda^0(\sigma)$, and in top degree $\mathcal P^-_r\Lambda^{\dim\sigma}(\sigma)=\mathcal P_{r-1}\Lambda^{\dim\sigma}(\sigma)$. We use the convention that a space of polynomial forms is zero if its polynomial degree is negative or its form degree lies outside $\{0,\dots,\dim \sigma\}$. For every simplex $\sigma$, we also set $\mathcal P_0^-\Lambda^0(\sigma)=\mathcal P_0\Lambda^0(\sigma) =\operatorname{span}\{1\}$ and $\mathcal P_0^-\Lambda^k(\sigma)=0$ for $k>0$. The lowest order trimmed spaces $\mathcal P_1^-\Lambda^k(\sigma)$ are exactly the Whitney forms \cite{Whitney:57} on a simplex.  
We record some well-known facts about these two families (see e.g. \cite[Sections 3-4]{Arnold-Falk-Winther:06}) that will be used later.

\begin{proposition}[Stability under the exterior derivative]
\label{prop:prelim--extder-polynomial}
The exterior derivative maps $\mathcal P_r\Lambda^k(\sigma)$ to $\mathcal P_{r-1}\Lambda^{k+1}(\sigma)$, and $\mathcal P_r^-\Lambda^k(\sigma)$ to $\mathcal P_r^-\Lambda^{k+1}(\sigma)$.
\end{proposition}

\begin{proposition}[Stability under the trace]
\label{prop:prelim--trace-polynomial}
Let $\tau\subset\sigma$ be a flat simplex (not necessarily a face of $\sigma$). Then $\tr_\tau$ maps $\mathcal P_r\Lambda^k(\sigma)$ to $\mathcal P_r\Lambda^k(\tau)$, and $\mathcal P_r^-\Lambda^k(\sigma)$ to $\mathcal P_r^-\Lambda^k(\tau)$.
\end{proposition}

We shall also need the subspaces with vanishing trace. For an $\ell$-simplex $\sigma$ and $0\leq k\leq\ell-1$, set
\begin{equation}
    \mathring{\mathcal P}_r\Lambda^k(\sigma)\coloneqq\set{u\in\mathcal P_r\Lambda^k(\sigma):\tr_{\partial\sigma}u=0}, \qquad
    \mathring{\mathcal P}_r^-\Lambda^k(\sigma)\coloneqq\set{u\in\mathcal P_r^-\Lambda^k(\sigma):\tr_{\partial\sigma}u=0}.
\end{equation}
For $k=\ell$ it holds $\mathring{\mathcal P}_r\Lambda^\ell(\sigma) = \mathcal P_r\Lambda^\ell(\sigma)$ and $\mathring{\mathcal P}_r^-\Lambda^\ell(\sigma) = \mathcal P_r^-\Lambda^\ell(\sigma) = \mathcal P_{r-1}\Lambda^\ell(\sigma)$.

\subsection{Global spaces, boundary conditions, and degrees of freedom}
\label{section:prelim--global-spaces}

The global spaces are assembled by requiring single-valued traces across faces. For any $u \in \bigoplus_{\sigma\in \mathcal T_h^n}\mathcal P_r\Lambda^k(\sigma)$, the discontinuous space of piecewise $\mathcal P_r\Lambda^k$-forms, we write $u_\sigma$ for the component of $u$ on $\sigma$, and we set
\begin{equation}
\label{eq:prelim--global-assembly}
\mathcal P_r\Lambda^k(\mathcal T_h) \coloneqq \set[\bigg]{u\in \bigoplus_{\sigma\in \mathcal T_h^n}\mathcal P_r\Lambda^k(\sigma)\ :\ \tr_{\sigma\cap\sigma'} u_\sigma = \tr_{\sigma\cap \sigma'} u_{\sigma'} \quad \forall \sigma,\sigma'\in\mathcal T_h^n },
\end{equation}
the trace conditions being vacuous when $\sigma\cap\sigma'=\emptyset$. The trimmed space $\mathcal P_r^-\Lambda^k(\mathcal T_h)$ is defined similarly. By \cref{prop:prelim--trace-polynomial}, the value $\tr_\tau u$ is then unambiguously defined for every $u\in \mathcal P_r\Lambda^k(\mathcal T_h)$ and $\tau\in\mathcal T_h$. 
Single-valued traces are precisely the condition for a piecewise smooth form to be globally $H\Lambda^k$-regular, so that \eqref{eq:prelim--global-assembly} may equivalently be written as a conforming subspace \cite[Section 5]{Arnold-Falk-Winther:06}:
\begin{equation}
\label{eq:prelim--conformity}
\mathcal P_r\Lambda^k(\mathcal T_h) = \set{ u\in H\Lambda^k(\Omega) : \tr_\sigma u\in \mathcal P_r\Lambda^k(\sigma) \quad \forall \sigma\in\mathcal T_h^n },
\end{equation}
and analogously for $\mathcal P_r^-\Lambda^k(\mathcal T_h)$. The global spaces with homogeneous boundary conditions are
\begin{equation}
\label{eq:prelim--bc-spaces}
\mathring{\mathcal P}_r\Lambda^k(\mathcal T_h)
\coloneqq\set{u\in\mathcal P_r\Lambda^k(\mathcal T_h):\tr_{\partial\Omega}u=0}, \qquad
\mathring{\mathcal P}_r^-\Lambda^k(\mathcal T_h)
\coloneqq\set{u\in\mathcal P_r^-\Lambda^k(\mathcal T_h):\tr_{\partial\Omega}u=0}.
\end{equation}
For $k=n$ it holds $\mathring{\mathcal P}_r\Lambda^n(\mathcal T_h) = \mathcal P_r\Lambda^n(\mathcal T_h)$ and $\mathring{\mathcal P}_r^-\Lambda^n(\mathcal T_h) = \mathcal P_r^-\Lambda^n(\mathcal T_h) =  \mathcal P_{r-1}\Lambda^n(\mathcal T_h)$.

To work with these spaces computationally, one uses degrees of freedom (dof). The relevant ones here are the canonical dof of \cite{Hiptmair:99, Arnold-Falk-Winther:06}, which pair a $k$-form against polynomial forms on the faces of complementary degree. 

\begin{proposition}[Canonical dof for the trimmed spaces]
\label{prop:prelim--dof-trimmed-polynomial}
Let $r\geq1$. For $u\in\mathcal P_r^-\Lambda^k(\mathcal T_h)$, the set of dof given by
\begin{equation}
\label{eq:prelim--dof-trimmed-polynomial}
u\mapsto \int_\tau \tr_\tau u\wedge q \qquad \forall q\in \mathcal P_{r+k-\ell-1}\Lambda^{\ell-k}(\tau),
\end{equation}
for all $\tau \in \mathcal T_h^\ell$ and all $\ell\geq k$, is unisolvent. The dof for the Whitney forms $u\in\mathcal P_1^-\Lambda^k(\mathcal T_h)$ reduce to integration over the $k$-faces,
\begin{equation}
u\mapsto \int_\tau \tr_\tau u \qquad \forall \tau \in \mathcal T_h^k.
\end{equation}
\end{proposition}

\begin{proposition}[Canonical dof for the full spaces]
\label{prop:prelim--dof-full-polynomial}
Let $r\geq1$. For $u\in\mathcal P_r\Lambda^k(\mathcal T_h)$, the set of dof given by
\begin{equation}
u\mapsto \int_\tau \tr_\tau u\wedge q \qquad \forall q\in \mathcal P_{r+k-\ell}^-\Lambda^{\ell-k}(\tau),
\end{equation}
for all $\tau \in \mathcal T_h^\ell$ and all $\ell\geq k$, is unisolvent.
\end{proposition}

Both statements are proved on a single simplex in \cite[Section 4]{Arnold-Falk-Winther:06}; see also \cite[Proposition 3.7]{christiansen-rapetti:15}. The global versions follow because the dof attached to a face $\tau$ depend on $u$ only through $\tr_\tau u$, which is single valued by \eqref{eq:prelim--global-assembly}. We note that the families of test spaces $\bigl(\mathcal P_{r+k-\ell-1}\Lambda^{\ell-k}(\tau)\bigr)_{k=0}^n$ and $\bigl(\mathcal P^-_{r-\ell}\Lambda^{\ell-k}(\tau)\bigr)_{k=0}^n$, of the dof for the trimmed and full complexes respectively, are themselves stable under the trace and the exterior derivative, by \cref{prop:prelim--extder-polynomial,prop:prelim--trace-polynomial}.

\subsection{Cohomology and homotopy/Poincar\'e operators}
\label{section:prelim--cohomology}

By \cref{prop:prelim--extder-polynomial}, both families of global spaces are stable under the exterior derivative and so we have the complexes:
\begin{subequations}
\label{eq:prelim--complexes}
\begin{align}
  &\begin{tikzcd}[column sep = large, ampersand replacement=\&]
    \mathcal P_r\Lambda^0(\mathcal T_h) \arrow[r, "d"] \&
    \mathcal P_{r-1}\Lambda^1(\mathcal T_h) \arrow[r, "d"] \&
    \cdots \arrow[r, "d"] \&
    \mathcal P_{r-n}\Lambda^n(\mathcal T_h),
  \end{tikzcd} \label{eq:prelim--complexes-a} \\
  &\begin{tikzcd}[column sep = large, ampersand replacement=\&]
    \mathcal P_r^-\Lambda^0(\mathcal T_h) \arrow[r, "d"] \&
    \mathcal P_r^-\Lambda^1(\mathcal T_h) \arrow[r, "d"] \&
    \cdots \arrow[r, "d"] \&
    \mathcal P_r^-\Lambda^n(\mathcal T_h).
  \end{tikzcd} \label{eq:prelim--complexes-b}
\end{align}
\end{subequations}
These are examples of what are known as \emph{Hilbert complexes} \cite{Arnold:18}. 

For our purposes, a Hilbert complex $(V^\bullet, d)$ is a sequence of finite dimensional vector spaces $V^0,\cdots, V^n$, each equipped with an inner product, and connected by linear maps $d^k:V^k\to V^{k+1}$ satisfying $d^{k+1}\circ d^k=0$. It is common to omit the superscript in $d^k$ when the space it is acting on is clear from context. 

\begin{definition}[Cohomology projection]
\label{def:prelim--cohomology-projection}
Let $(V^\bullet,d)$ be a Hilbert complex. A \emph{cohomology projection} is a graded linear operator $\pi^k:V^k \to V^k$ such that 
\begin{itemize}
    \item[i)] $\pi^k \circ \pi^k = \pi^k$;
    \item[ii)] $d^k\pi^k = \pi^{k+1} d^k = 0$;
    \item[iii)] If $z\in \ker d^k$, then $z-\pi^kz\in \im d^{k-1}$.
\end{itemize}
\end{definition}

The role of $\pi^\bullet$ is to encode cohomology. The subcomplex $\ker\pi^\bullet$ has trivial cohomology, while $\im\pi^k\subset \ker d$ maps isomorphically onto the $k$'th cohomology group via the quotient map
\begin{equation}
    \im \pi^k \cong  \ker d^k/ \im d^{k-1}. 
\end{equation}
In other words, the operator $\pi^k$ selects a representative of the $k$'th cohomology group.

\begin{definition}[Homotopy operator]
\label{def:prelim--homotopy-operator}
Let $(V^\bullet,d)$ be be a Hilbert complex, and $\pi^k$ a cohomology projection. A \emph{homotopy operator} on $V^\bullet$ associated with $\pi^\bullet$ is a graded map $\mathscr H^k:V^k\to V^{k-1}$ satisfying the homotopy identity
\begin{equation}
\label{eq:prelim--homotopy}
d^{k-1}\mathscr H^ku +\mathscr H^{k+1}d^ku = u - \pi^k u \qquad \forall u\in V^k.
\end{equation}
Here we always use the convention that $\mathscr H^0=\mathscr H^{n+1}=0$.
\end{definition}

\begin{definition}[Poincar\'e operator]
\label{def:prelim--poincare}
    In the context of differential forms, a Poincar\'e operator $\mathscr P^\bullet$ is a homotopy operator associated with a cohomology projection $\pi^\bullet$ such that $\pi^k=0$ for all $k>0$ and the range of $\pi^0$ is the space of constant $0$-forms.
\end{definition}

The classical example is the Poincar\'e operator of the smooth de Rham complex on a domain $\Omega$ which is star-shaped with respect to a point $p_0$. For a smooth $k$-form $u$, $k\geq1$, it is given by
\begin{equation}
\label{eq:prelim--classical-poincare}
(\mathscr P^k_{p_0}u)_x(v_1,\dots,v_{k-1})\coloneqq\int_0^1t^{k-1}\,u_{p_0+t(x-p_0)}\bigl(x-p_0,v_1,\dots,v_{k-1}\bigr)\,dt,
\end{equation}
and it satisfies \eqref{eq:intro--homotopy} with $\pi^0u=u(p_0)$ \cite[Lemma 17.9]{Lee:13}.


\section{Construction by solving global discrete \texorpdfstring{$\mathfrak{B}^*$}{B*} problems}
\label{section:global-construction}

In this section, we show how a discrete Poincar\'e operator arises from solving $\mathfrak B^*$ problems (see e.g. \cite[Section 3.2.1]{Arnold-Falk-Winther:10}) over the global discrete spaces. We begin by formulating the $\mathfrak B^*$ problem for an abstract complex.

\subsection{Homotopy operators from \texorpdfstring{$\mathfrak B^*$}{B*} problems}
\label{section:prelim--abstract-bstar}

Let $(V^\bullet,d)$ be a Hilbert complex, with $(\cdot,\cdot)$ denoting the inner product on $V^k$. We introduce the spaces
\begin{subequations}
\label{eq:prelim--abstract-subspaces}
\begin{alignat}{2}
\mathfrak ZV^k&\coloneqq\set{u\in V^k: du=0},\\
\mathfrak HV^k&\coloneqq\set{z\in\mathfrak ZV^k: (z,dv)=0\quad\forall v\in V^{k-1}}.
\end{alignat}
\end{subequations}

\begin{definition}
\label{def:prelim--abstract-bstar}
Given $u\in V^k$ with $k\geq 1$, the \emph{$\mathfrak B^*$ problem} for $u$ is: find $v\in V^{k-1}$ such that
\begin{subequations}
\label{eq:prelim--abstract-bstar}
\begin{alignat}{2}
(dv,dv')&=(u,dv')\qquad&&\forall v'\in V^{k-1},\label{eq:prelim--abstract-bstar-a}\\
(v,z)&=0\quad&&\forall z\in\mathfrak ZV^{k-1}.\label{eq:prelim--abstract-bstar-b}
\end{alignat}
The first condition \eqref{eq:prelim--abstract-bstar-a} uniquely determines $dv$ as the orthogonal projection of $u$ onto $dV^{k-1}$. The second condition \eqref{eq:prelim--abstract-bstar-b} fixes the choice of representation $v$ in the pre-image of $d$. In practice, to find $v$ one would solve a mixed formulation using the decomposition $\mathfrak ZV^k=\mathfrak{H} V^k\oplus \im d^{k-1}$, rather than trying to work with the space $\mathfrak Z V^{k-1}$ directly (see \cite[Equation (4.16) and Remark 4.10]{Arnold:18}).
\end{subequations}
\end{definition}

\begin{proposition}[Homotopy operator from the $\mathfrak B^*$ problem]
\label{prop:prelim--abstract-bstar}
Let $(V^\bullet,d)$ be a Hilbert complex. For $1\leq k\leq n$ let $\mathscr H^k_V:V^k\to V^{k-1}$ be the map assigning to $u\in V^k$ the solution $v\eqqcolon\mathscr H^k_V u$ of the $\mathfrak B^*$ problem. Let $\pi^k_V:V^k\to V^k$ be the orthogonal projection onto $\mathfrak HV^k$. Then
\begin{equation}
\label{eq:prelim--abstract-homotopy}
d\mathscr H^k_Vu +\mathscr H^{k+1}_Vdu=u-\pi^k_Vu\qquad\forall u\in V^k,\ 0\leq k\leq n,
\end{equation}
and $\pi^\bullet_V$ is a cohomology projection.

\end{proposition}
\begin{proof}
    Using the Hodge decomposition $\mathfrak ZV^k=\mathfrak HV^k\oplus \im d^{k-1}$, it is straightforward to check that $\pi^k_V$ is a cohomology projection. We show \eqref{eq:prelim--abstract-homotopy} by checking that $e\coloneqq u-\pi^k_Vu - d\mathscr H^k_Vu -\mathscr H^{k+1}_V du$ is zero. Plugging $de=du - d\mathscr H^{k+1}_V du$ into \eqref{eq:prelim--abstract-bstar-a} yields $de=0$, so it suffices to show that $e\perp \mathfrak ZV^k$. 
    
    For any $z\in \mathfrak ZV^k$ we may write $z=h+dq$, with $h\in \mathfrak HV^k$ and $q\in V^{k-1}$. Testing against $h$ yields
    \begin{equation}
        (e,h) = (u-\pi^k_Vu, h) - (d\mathscr H^k_Vu, h) - (\mathscr H^{k+1}_Vdu, h).
    \end{equation}
    The first term is zero since $\pi^k_Vu$ is the orthogonal projection of $u$ onto $\mathfrak HV^k$, the second is zero since $h\in\mathfrak HV^k$, and the third term is zero by \eqref{eq:prelim--abstract-bstar-b}. Similarly, testing against $dq$ we get
    \begin{equation}
        (e,dq) = (u-d\mathscr H^k_Vu, dq) -  (\pi^k_Vu, dq) - (\mathscr H^{k+1}_Vdu, dq).
    \end{equation}
    The first term is zero by \eqref{eq:prelim--abstract-bstar-a}, the second since $\pi^k_Vu\in \mathfrak H V^k$, and the third term is zero by \eqref{eq:prelim--abstract-bstar-b}.
\end{proof}

We apply \cref{def:prelim--abstract-bstar} to obtain a simple discrete Poincar\'e operator for the trimmed polynomial forms. Given $u\in \mathcal P_r^-\Lambda^k(\mathcal T_h)$, let $\mathscr H^k_ru\in\mathcal P_r^-\Lambda^{k-1}(\mathcal T_h)$ be the unique solution of the discrete $\mathfrak B^*$ problem
\begin{subequations}
\label{eq:prelim--hodge-laplace}
\begin{alignat}{2}
(d\mathscr H^k_ru, dv')_\Omega &= (u, dv')_\Omega \qquad &&\forall v'\in\mathcal P_r^-\Lambda^{k-1}(\mathcal T_h), \\
(\mathscr H^k_r u, z)_\Omega  &= 0 &&\forall z\in \mathfrak Z \mathcal P_r^-\Lambda^{k-1}(\mathcal T_h).
\end{alignat}
By \cref{prop:prelim--abstract-bstar} this defines a homotopy operator for the complex $V^\bullet=\mathcal P_r^-\Lambda^\bullet(\mathcal T_h)$. In particular, if $\Omega$ is contractible then $\mathscr H_r^\bullet$ is a Poincar\'e operator. 
\end{subequations}

Due to its simplicity, \eqref{eq:prelim--hodge-laplace} is an appropriate operator for most applications. It is, however, a global problem posed over the entire discrete space, and it is far from the only construction. In the following section, we will show how a decomposition allows one to obtain a more numerically favourable construction.


\section{Construction from Whitney-remainder decomposition}
\label{section:local-construction}

We will decompose the trimmed polynomial complex \eqref{eq:prelim--complexes-b} into its Whitney subcomplex and a complementary complex. This idea was used in \cite[Section 5.2]{DiPietro-Droniou-Pitassi:25} to develop bounded discrete Poincar\'e operators for analysis purposes. We briefly review their construction, and introduce an alternative way to handle the complementary complex via discrete $\mathfrak{B}^*$ problems. 

We note that the idea of splitting the higher order complex in this way has also been used in the development of bounded commuting projections \cite{falk-winther:14, arnold-guzman:21} where the higher order complement involves only local problems. While we present the construction only for the trimmed polynomial complex, the idea of the construction extends to the full polynomial complex \eqref{eq:prelim--complexes-a} with appropriate modification.

\subsection{Decomposition of \texorpdfstring{$\mathcal P_r^-\Lambda^k(\mathcal T_h)$}{the trimmed polynomial forms}}
\label{section:lc--decomposition}

For $u\in\mathcal P_r^-\Lambda^k(\mathcal T_h)$, the canonical interpolant onto the Whitney forms $I^k u\in\mathcal P_1^-\Lambda^k(\mathcal T_h)$ is defined by
\begin{equation}
\label{eq:lc--whitney-interpolant}
\int_\tau\tr_\tau I^k u=\int_\tau\tr_\tau u
\qquad\forall\tau\in\mathcal T_h^k.
\end{equation}
The same conditions on the $k$-faces of a simplex $\sigma$ yield local interpolators $I_\sigma^k:\mathcal P_r^-\Lambda^k(\sigma)\to \mathcal P_1^-\Lambda^k(\sigma)$. They are related by $I_\sigma^k \tr_\sigma u = \tr_\sigma I^k u$. We note the following well-known properties of $I^k$ and $I_\sigma^k$ (see \cite[Theorem 5.2]{Arnold-Falk-Winther:06}).

\begin{proposition}[Properties of the canonical interpolant]
\label{lem:lc--interpolant}
It holds,\\
\textup{(i)} $I^k$ and $I_\sigma^k$ are projections: $I^k \circ I^k = I^k$,\hspace{.5em} $I_\sigma^k \circ I_\sigma^k = I_\sigma^k$.\\
\textup{(ii)} For any $u\in\mathcal P_r^-\Lambda^k(\mathcal T_h)$, respectively $u\in \mathcal P_r^-\Lambda^k(\sigma)$,
\begin{equation}
\label{eq:lc--interpolant-commutes}
dI^k u =I^{k+1}d u, \quad dI_\sigma^k u =I_\sigma^{k+1}d u.
\end{equation}
\textup{(iii)} For all $\sigma\in\mathcal T_h$ and $\tau\triangleleft\sigma$, 
\begin{equation}
\label{eq:lc--interpolant-trace}
\tr_\tau I_\sigma^k u = I_\tau^k\tr_\tau u.
\end{equation}
\textup{(iv)} For all $\tau\in \partial \mathcal T_h$ and $\sigma \triangleright \tau$, if $\tr_\tau u = 0$, then $\tr_\tau I^k u = \tr_\tau I_\sigma^k u = 0$.
\end{proposition}

\begin{definition}[Whitney remainder]
    
The global and local remainder spaces are given by
\begin{subequations}
\label{eq:lc--remainder-spaces}
\begin{align}
\mathcal M_r^-\Lambda^k(\mathcal T_h)
&\coloneqq\ker I^k
=\set[\bigg]{u\in\mathcal P_r^-\Lambda^k(\mathcal T_h):
\int_\tau\tr_\tau u=0\quad\forall\tau\in\mathcal T_h^k},
\label{eq:lc--remainder-global}\\
\mathcal M_r^-\Lambda^k(\sigma)
&\coloneqq\ker I_\sigma^k
=\set[\bigg]{u\in\mathcal P_r^-\Lambda^k(\sigma):
\int_\tau\tr_\tau u=0\quad
\forall\tau\triangleleft\sigma,\, \dim\tau=k}.
\label{eq:lc--remainder-local}
\end{align}
\end{subequations}
We set the spaces to zero outside the admissible form degrees. The versions with homogeneous boundary conditions, $\mathring{\mathcal M}_r^-\Lambda^k(\mathcal T_h)$ and $\mathring{\mathcal M}_r^-\Lambda^k(\sigma)$, are defined as before.
\end{definition}

By \eqref{eq:lc--interpolant-commutes}, the remainder spaces form subcomplexes of their corresponding polynomial complexes. Moreover, they have been constructed so that the cohomology of $\mathcal M_r^-\Lambda^k(\sigma)$, as well as that of $\mathring{\mathcal M}_r^-\Lambda^k(\sigma)$, vanishes in every degree. That is, the following sequences are exact 
\begin{subequations}
\begin{equation}
\label{eq:lc--local-exact}
0\longrightarrow{\mathcal M}_r^-\Lambda^0(\sigma)
\xrightarrow d\cdots\xrightarrow d
{\mathcal M}_r^-\Lambda^\ell(\sigma)
\longrightarrow0;
\end{equation}
\begin{equation}
\label{eq:lc--local-exact-bubble}
0\longrightarrow\mathring{\mathcal M}_r^-(\sigma)
\xrightarrow d\cdots\xrightarrow d
\mathring{\mathcal M}_r^-\Lambda^\ell(\sigma)
\longrightarrow0;
\end{equation}
\end{subequations}
where $\ell=\dim \sigma$. Indeed, since $\id - I^k_\sigma$ is a cochain projection, the cohomology spaces of \eqref{eq:lc--local-exact} and \eqref{eq:lc--local-exact-bubble} are at most as large as the corresponding de Rham cohomology spaces on $\sigma$ \cite[Proposition 2.1]{Arnold:18}. As for the potentially non-zero entries, \eqref{eq:lc--local-exact} is exact at degree $k=0$ since $\mathcal M_r^-\Lambda^0(\sigma)$ does not contain constants, and \eqref{eq:lc--local-exact-bubble} is exact at degree $k=\ell$ since elements of $\mathring{\mathcal M}_r^-\Lambda^{\ell}(\sigma)$ have zero average.

We conclude this subsection by noting that the following are direct sums of complexes:
\begin{align}
\mathcal P_r^-\Lambda^\bullet(\mathcal T_h)
&=\mathcal P_1^-\Lambda^\bullet(\mathcal T_h)
\oplus\mathcal M_r^-\Lambda^\bullet(\mathcal T_h),
\label{eq:lc--wr-decomposition-global}\\
\mathring{\mathcal P}_r^-\Lambda^\bullet(\mathcal T_h)
&=\mathring{\mathcal P}_1^-\Lambda^\bullet(\mathcal T_h)
\oplus\mathring{\mathcal M}_r^-\Lambda^\bullet(\mathcal T_h).
\label{eq:lc--relative-decomposition}
\end{align}
The Whitney component of $u$ is $I^ku$ and the remainder component is $(\id-I^k)u$ for both decompositions.

\subsection{Main idea of the decomposition Poincar\'e operator}

\begin{proposition}[Abstract construction]
\label{prop:lc--abstract-construction}
Let $\mathscr H_{\mathrm W}^k:\mathcal P_1^-\Lambda^k(\mathcal T_h) \longrightarrow\mathcal P_1^-\Lambda^{k-1}(\mathcal T_h)$ be a homotopy operator on the Whitney complex, with cohomology projection $\pi_{\mathrm W}^\bullet$. Suppose that $B^k:\mathcal M_r^-\Lambda^k(\mathcal T_h) \longrightarrow\mathcal M_r^-\Lambda^{k-1}(\mathcal T_h)$ satisfies
\begin{equation}
\label{eq:lc--ag-homotopy}
dB^k u +B^{k+1}d u = u
\quad\forall u \in \mathcal M_r^-\Lambda^k(\mathcal T_h).
\end{equation}
Then
\begin{equation}
\label{eq:lc--higher-order-homotopy}
\mathscr H_{r,\mathrm{dec}}^ku
\coloneqq\mathscr H_{\mathrm W}^k I^ku+B^k(\id-I^k)u
\end{equation}
is a homotopy operator on $\mathcal P_r^-\Lambda^\bullet(\mathcal T_h)$ with cohomology projection
$\pi_{r,\mathrm{dec}}^k\coloneqq\pi_{\mathrm W}^k I^k$.
\end{proposition}

\begin{proof}
Using \eqref{eq:lc--interpolant-commutes}, the homotopy identity for $\mathscr H_{\mathrm W}^k$, and \eqref{eq:lc--ag-homotopy} gives
\begin{align}
d\mathscr H_{r,\mathrm{dec}}^k+\mathscr H_{r,\mathrm{dec}}^{k+1}d
&=(d\mathscr H_{\mathrm W}^k+\mathscr H_{\mathrm W}^{k+1}d)I^k
 +(dB^k+B^{k+1}d)(\id-I^k)\\ \nonumber
&=(\id-\pi_{\mathrm W}^k)I^k+\id-I^k
=\id-\pi_{\mathrm W}^kI^k.
\end{align}
This identity also shows that $\pi_{r,\mathrm{dec}}^\bullet$ represents the cohomology of the higher order complex. Since $I^k$ is the identity on the Whitney space, the remaining properties of a cohomology projection follow directly from those of $\pi_{\mathrm W}^\bullet$.
\end{proof}

If $\mathscr H_{\mathrm W}^\bullet$ is a Poincar\'e operator, then so is $\mathscr H_{r,\mathrm{dec}}^\bullet$, in which case we may write $\mathscr P_{r,\mathrm{dec}}^\bullet$ in place of $\mathscr H_{r,\mathrm{dec}}^\bullet$. On a general domain, the $r=1$ instance of \eqref{eq:prelim--hodge-laplace} supplies a possible choice of $\mathscr H_{\mathrm W}^\bullet$ and $\pi_{\mathrm W}^\bullet$. Other choices include those discussed in \cite{guzman-potu:26}. Hence, we now focus on construction of the operator $B^k$.

\subsection{Di Pietro--Droniou--Pitassi construction of \texorpdfstring{$B^k$}{Bk}}
\label{section:lc--DDP-construction}

We now briefly detail how $B^k$ is constructed in \cite[Section 5]{DiPietro-Droniou-Pitassi:25}. The essential idea is: given $u\in \mathcal M_r^-\Lambda^k(\mathcal T_h)$, we define $B^k u \in \mathcal M_r^-\Lambda^{k-1}(\mathcal T_h)$ by setting the dof of \Cref{prop:prelim--dof-trimmed-polynomial} so that the homotopy identity \eqref{eq:lc--ag-homotopy} holds. This is done by induction as follows. 

Set $B^0 = B^{n+1} = 0$ and let $k\in\{1,\dots,n\}$. Suppose that $B^{k+1}$ satisfying \eqref{eq:lc--ag-homotopy} has already been defined.  Let $u\in \mathcal M_r^-\Lambda^k(\mathcal T_h)$. We set the dof of $B^k$, as given in \eqref{eq:prelim--dof-trimmed-polynomial}, by iterating over the dimension of $\tau \in \mathcal T_h^\ell$. First, for $\tau \in \mathcal T_h^{k-1}$ we set 
\begin{equation}
    \int_{\tau} (B^k u)\wedge q \coloneqq 0, \quad \forall q\in \mathcal P_{r-1}\Lambda^0(\tau). 
\end{equation}
Next, suppose the dofs corresponding to $\tau'\in \mathcal T_h^j$ for all $j\in\{k-1, \dots, \ell-1\}$ have been set. In particular, this means that $\tr_{\partial \tau} B^k u$ has been determined for all $\tau\in \mathcal T_h^\ell$. Note the decomposition $\mathcal P_t \Lambda^s(\sigma) = d \kappa_a \mathcal P_t \Lambda^s(\sigma) \oplus \kappa_a \mathcal P_{t-1}\Lambda^{s+1}(\sigma)$, where $\kappa_a$ is the Koszul operator on $\sigma$. This means that any $q\in \mathcal P_t \Lambda^s(\sigma)$ can be written $q = d\omega_q + \zeta_q$ with $\omega_q\in \kappa_a \mathcal P_t \Lambda^s(\sigma)$. For $\tau\in\mathcal T_h^\ell$ we set
\begin{equation}
    \int_\tau B^k u \wedge q = (-1)^k\int_{\tau} (u- B^{k+1}du) \wedge \omega_q + (-1)^{k-1}\int_{\partial\tau} B^k u \wedge \omega_q, \quad \forall q\in \mathcal P_{r+k-\ell-2}\Lambda^{\ell-k+1}(\tau).
\end{equation}
Iterating this process of setting dof over all $\ell$ fully defines $B^k u$. One can then show that the homotopy identity \eqref{eq:lc--ag-homotopy} holds by considering $e:= u - d B^{k} u - B^{k+1} du$ and showing that all the dof of \eqref{eq:prelim--dof-trimmed-polynomial} vanish for $e$. We omit this proof for brevity and point to \cite[Section 5]{DiPietro-Droniou-Pitassi:25} for details.

\subsection{Discrete \texorpdfstring{$\mathfrak{B}^*$}{B*} problem construction of \texorpdfstring{$B^k$}{Bk}}
\label{section:lc--Bk-construction}

We provide an alternative construction of $B^k$ based on solving \eqref{eq:prelim--hodge-laplace} over the local spaces. This perspective has the advantage of being easy to implement in standard software (see \cref{rem:lc--cost} and \cref{section:implementation}). For each $\sigma\in\mathcal T_h^m$ we will construct a local operator on $B_\sigma^k$ satisfying the homotopy identity. We use induction to ensure that everything glues together nicely, so that later we may combine these local operators to form one global operator $B^k$. 

For $\dim \sigma = 0$, as well as for $\dim \sigma \leq k-1$, we set $B^k_\sigma:\mathcal M_r^-\Lambda^k(\sigma)\to\mathcal M_r^-\Lambda^{k-1}(\sigma)$ to be the zero map. Otherwise, assuming the maps $B^k_\tau$ have already been defined on all facets $\tau\triangleleft \sigma$ (all co-dimension 1 faces of $\sigma$), we define $B^k_\sigma u\in\mathcal M_r^-\Lambda^{k-1}(\sigma)$ to be the solution to the problem
\begin{subequations}
\label{eq:lc--local-problem}
\begin{alignat}{2}
(dB^k_\sigma u,dv')_\sigma + (d\lambda, v')_\sigma 
&=(u,dv')_\sigma\quad&&\forall v'\in \mathring{\mathcal M}_r^-\Lambda^{k-1}(\sigma),\label{eq:lc--local-problem-a}
\\
(B^k_\sigma u,dq)_\sigma - (\lambda, q)_\sigma
&=0\quad&&\forall q\in\mathring{\mathcal M}_r^-\Lambda^{k-2}(\sigma),\label{eq:lc--local-problem-b}
\\
\tr_\tau B^k_\sigma u
&= B^k_\tau \tr_\tau u; \quad&&\text{$\forall$ facets $\tau\triangleleft \sigma$.} \label{eq:lc--local-problem-c}
\end{alignat}
\end{subequations}
Here $\lambda\in \mathring{\mathcal M}_r^-\Lambda^{k-2}(\sigma)$ is a free auxiliary variable, giving us a mixed formulation. The final solution will satisfy $\lambda=0$, so that \eqref{eq:lc--local-problem} reduces to a variant of the $\mathfrak B^*$ problem, but which is more amenable to computation than $\eqref{eq:prelim--abstract-bstar}$ \cite{ern-guermond:21-2,Ham-et-al:23}. Note that we impose boundary conditions on the test spaces in \eqref{eq:lc--local-problem-a}-\eqref{eq:lc--local-problem-b}.

\begin{proposition}
\label{prop:lc--local-operators}
The problem \eqref{eq:lc--local-problem} gives a square linear system which admits a unique solution, and the following holds for all $\tau \triangleleft \sigma$ (not just for facets) 
\begin{equation}
\label{eq:lc--local-operators-trace}
    \tr_\tau B^k_\sigma u = B^k_\tau\tr_\tau u.
\end{equation}
\end{proposition}

\begin{proof}
Both statements are trivial for $\dim \sigma=0$. We proceed by induction, demonstrating first the well posedness of \eqref{eq:lc--local-problem}, and then showing that \eqref{eq:lc--local-operators-trace} follows from the boundary condition.

Let $\sigma\in\mathcal T^\ell_h$, and suppose the statement holds for all $\tau\in\mathcal T_h$ with $\dim \tau \leq \ell-1$. For any pair of facets $\tau,\tau'\triangleleft \sigma$ the induction hypothesis \eqref{eq:lc--local-operators-trace} yields
\begin{equation}
\tr_{\tau\cap \tau'}B^k_\tau\tr_\tau u = \tr_{\tau\cap \tau'} B^k_{\tau'} \tr_{\tau'} u.
\end{equation}
Thus $g_\tau=B^k_\tau u$ defines a form $g\in \mathcal M_r^-\Lambda^{k-1}(\partial \sigma)$. By zeroing the internal dof of \cref{prop:prelim--dof-trimmed-polynomial} we can extend $g$ to some $\tilde g\in \mathcal M_r^-\Lambda^{k-1}(\sigma)$. Let $v\in \mathring{\mathcal M}_{r}^-\Lambda^{k-1}(\sigma)$ be the solution to the $\mathfrak B^*$ problem for $u-d\tilde g$, as posed over $\mathring{\mathcal M}_r^-\Lambda^\bullet$, and let $\Pi\tilde g$ be the $L^2$ projection of $\tilde g$ onto $d\mathring{\mathcal M}_r^-\Lambda^{k-2}(\sigma)$. Then $B^k_\sigma u \coloneqq  v+\tilde g-\Pi \tilde g$ satisfies \eqref{eq:lc--local-problem} with $\lambda=0$. That this is the only such solution is due to the well-posedness of the Hodge Laplace problem \cite{Arnold:18}, using that $\mathfrak H\mathring{\mathcal M}_r^-\Lambda^{k-1}(\sigma)= \set{0}$.

Lastly we check that \eqref{eq:lc--local-operators-trace} holds. Given $\tau \triangleleft \sigma$ with $\tau\ne \sigma$, there exists a facet $F\triangleleft \sigma$ such that $\tau\triangleleft F\triangleleft \sigma$. Using \eqref{eq:lc--local-problem-c} together with the induction hypothesis we compute
\begin{equation}
    \tr_\tau B^k_\sigma u = \tr_{\tau}\tr_F B^k_\sigma u = B^k_\tau \tr_\tau u.
\end{equation}
This concludes the proof.
\end{proof}

For $u\in\mathcal M_r^-\Lambda^k(\mathcal T_h)$, we define the global $B^ku$ cellwise by
\begin{equation}
\label{eq:lc--global-B}
\tr_\sigma B^ku\coloneqq B_\sigma ^k\tr_\sigma u
\qquad\forall \sigma\in\mathcal T_h^n.
\end{equation}
It follows from \eqref{eq:lc--local-operators-trace} that $B^ku$ has single-valued traces, and thus belongs to $\mathcal M_r^-\Lambda^{k-1}(\mathcal T_h)$. We note also that $B^k$ preserves homogeneous boundary conditions.

\begin{proposition}[Homotopy formula]
\label{prop:lc--global-contraction}
The operator $B^k:\mathcal M_r^-\Lambda^k(\mathcal T_h) \longrightarrow\mathcal M_r^-\Lambda^{k-1}(\mathcal T_h)$ given by \eqref{eq:lc--global-B}, satisfies the homotopy formula \eqref{eq:lc--ag-homotopy}
\begin{equation*}
    dB^k u + B^{k+1}du = u.
\end{equation*}
\end{proposition}

\begin{proof}

It suffices to prove the identity holds locally on every $\sigma \in \mathcal T_h$. We do so by induction on the dimension of $\sigma$.

For $\dim \sigma\leq k$, and in particular for $\dim\sigma=0$, the statement is trivial as then $B^k_\sigma v=0$. Suppose the statement holds locally for all $\tau \in\mathcal T_h$ with $\dim \tau \leq \dim\sigma -1$, and let $u\in \mathcal M_r^-\Lambda^k(\sigma)$. The strategy is to show that $e\coloneqq u - dB^k_\sigma u - B^{k+1}_\sigma du$ is zero. By the induction hypothesis it holds $\tr_\tau e = 0$ for all facets $\tau \triangleleft \sigma$, so that $e\in\mathring{\mathcal M}_r^-\Lambda^k(\sigma)$. From here the argument that $e=0$ is identical to that in the proof of \cref{prop:prelim--abstract-bstar}.
\end{proof}

\begin{remark}[The full polynomial family]
\label{rem:lc--full-family}
Assume $r\geq n$, so that \eqref{eq:prelim--complexes-a} contains the Whitney complex as a subcomplex. The interpolant \eqref{eq:lc--whitney-interpolant} then makes sense on $\mathcal P_{r-k}\Lambda^k(\mathcal T_h)$ and still commutes with $d$, so \eqref{eq:lc--wr-decomposition-global} and \cref{prop:lc--abstract-construction} carry over with $\mathcal M_r^-\Lambda^k$ replaced by $\ker (I^k:\mathcal P_{r-k}\Lambda^k \to \mathcal P_1^-\Lambda^k)$. 
\end{remark}

\begin{remark}[Implementation]
\label{rem:lc--cost}
One may replace the $L^2$ inner product on $\sigma$ by the one on the reference simplex $\widehat\Delta^{\dim\sigma}$, by pulling back via the affine map $F_\sigma$ sending vertices to vertices in increasing order. The local operator then becomes $B^k_\sigma=(F_\sigma^*)^{-1}\widehat B^k_mF_\sigma^*$ for $\sigma\in\mathcal T_h^m$, and the reference operators $\widehat B^k_m$ need only be computed once. The only global object in \eqref{eq:lc--higher-order-homotopy} is $\mathscr H_{\mathrm W}^\bullet$, which is independent of $r$. The latter can be produced combinatorially, without solving a linear system, by the constructions of \cite{guzman-potu:26}.
\end{remark}

We conclude this section by computing stability estimates for the resulting operators.

\begin{proposition}[Stability of the local operator]
\label{prop:lc--local-stability}
For $\sigma\in\mathcal T_h$, let $h_\sigma$ be its diameter and $\rho_\sigma$ be the diameter of its inscribed ball. The mesh size is  $h=h(\mathcal T_h)\coloneqq\max_{\sigma\in\mathcal T_h^n}h_\sigma$ and the shape-regularity constant is $\gamma=\gamma(\mathcal T_h)\coloneqq\max_{\sigma\in\mathcal T_h^n}h_\sigma/\rho_\sigma$. For the operator $B^k$ given by $\eqref{eq:lc--global-B}$, there is a constant $C$ depending only on $r$, $n$, and $\gamma$ such that
\begin{equation}
\label{eq:lc--local-stability}
    \normdelim{B^k u}_{L^2(\Omega)} + h\normdelim{dB^k u}_{L^2(\Omega)} \leq Ch\normdelim{u}_{L^2(\Omega)}.
\end{equation}

\begin{proof}
Let $\sigma\in\mathcal T_h$. For $\tau\triangleleft \sigma$, let $L_\tau\colon\mathcal M_r^-\Lambda^k(\partial\tau)\to\mathcal M_r^-\Lambda^k(\tau)$ denote the lifting $g\mapsto \tilde g$ as described in the proof of \cref{prop:lc--local-operators}. That means $L_\tau g$ is obtained by prescribing the canonical dof \eqref{eq:prelim--dof-trimmed-polynomial}, matching those of $g$ on the proper faces of $\tau$, and set to zero otherwise. Standard scaling arguments \cite[Lemmas 11.1, 11.7 and 12.8]{Ern-Guermond:21} give the following three inequalities, which hold for all $\tau\triangleleft \sigma$ with constants depending only on $n$, $r$, and $\gamma$. 

For $w\in \mathring{\mathcal M}_r^-\Lambda^k(\tau)$, the discrete Poincar\'e inequality
\begin{equation}
\label{eq:app--poincare}
w\perp \mathfrak Z\mathring{\mathcal M}_r^-\Lambda^k(\tau) \implies \|w\|_{L^2(\tau)}\leq C\,h_\sigma\,\|dw\|_{L^2(\tau)}.
\end{equation}
For $g\in \mathcal M_r^-\Lambda^k(\partial \tau)$, stability of the lifting operator
\begin{equation}
\label{eq:app--lifting}
\normdelim{L_\tau g}_{L^2(\tau)} + h_\sigma\normdelim{dL_\tau g}_{L^2(\tau)} \leq Ch^{1/2}_\sigma\normdelim{g}_{L^2(\partial\tau)}.
\end{equation}
For $v\in \mathcal M_r^-\Lambda^k(\tau)$, the discrete trace inequality
\begin{equation}
\label{eq:app--trace}
\|\tr_{\partial \tau} v\|_{L^2(\partial \tau)}\leq C\,h_\sigma^{-1/2}\,\|v\|_{L^2(\tau)}.
\end{equation}
From here we show that \eqref{eq:lc--local-stability} holds locally by an induction argument.

For each $\tau \triangleleft \sigma$ we want to show that 
\begin{equation}
\label{eq:app--claim}
\|B^k_\tau v\|_{L^2(\tau)}+h_\sigma\|dB^k_\tau v\|_{L^2(\tau)}\leq C\,h_\sigma\,\|v\|_{L^2(\tau)},
\end{equation}
which should hold for all $v\in\mathcal M_r^-\Lambda^k(\tau)$ with $C$ depending only on $n$, $r$, and $\gamma$. For $\dim\tau=0$ the statement is trivial, so suppose $\dim \tau \geq 1$ and that \eqref{eq:app--claim} holds on all faces of strictly lower dimension. We write $B^k_\tau v= w+ L_\tau g - \Pi L_\tau g$, where $g\coloneqq \tr_{\partial\tau} B^k_\tau v$, and $\Pi$ denotes the $L^2$-orthogonal projection onto $d\mathring{\mathcal M}_r^-\Lambda^{k-2}(\tau)$. Then $w=B^k_\tau v-L_\tau g+\Pi L_\tau g$ lies in $\mathring{\mathcal M}_r^-\Lambda^k(\tau)$, and by \eqref{eq:lc--local-problem-b} it holds $w\perp \mathfrak Z\mathring{\mathcal M}_r^-\Lambda^{k-1}(\tau)$ (recall that $\mathfrak H\mathring{\mathcal M}_r^-\Lambda^{k-1}(\tau) = \set{0}$). 

By \eqref{eq:lc--local-operators-trace} it holds $g_{\tau'} = B^k_{\tau'} \tr_{\tau'} v$, and together with the induction hypothesis this yields
\begin{align}
\label{eq:lc--lifting-bound}
    \normdelim{L_\tau g}_{L^2(\tau)} + h_\sigma\normdelim{dL_\tau g}_{L^2(\tau)}
    &\leq Ch_\sigma^{1/2}\normdelim{g}_{L^2(\partial\tau)} \\ \nonumber
    &\leq Ch_\sigma^{1/2}\sum_{\text{facets }\tau'\triangleleft \tau}\normdelim{B^k_{\tau'} \tr_{\tau'} v}_{L^2(\tau')} \\ \nonumber
    &\leq Ch_\sigma^{3/2}\sum_{\text{facets }\tau'\triangleleft \tau}\normdelim{\tr_{\tau'} v}_{L^2(\tau')} \\ \nonumber
    &\leq Ch_\sigma\normdelim{v}_{L^2(\tau)}.
\end{align}
The first line of \eqref{eq:lc--lifting-bound} is \eqref{eq:app--lifting}, and the last line is \eqref{eq:app--trace}. As for $w$, the fact that it is orthogonal to $\mathfrak Z\mathring{\mathcal M}_r^-\Lambda^{k-1}(\tau)$ lets us apply the discrete Poincar\'e inequality \eqref{eq:app--poincare} to get $\normdelim{w}_{L^2(\tau)} \leq Ch_\sigma\normdelim{dw}_{L^2(\tau)}$. Plugging $dw$ into \eqref{eq:lc--local-problem-a}, and using the fact that $\lambda=0$, yields
\begin{equation}
    (dw, dw)_\tau = (u, dw) - (dL_\tau g,dw),
\end{equation}
and consequently
\begin{equation}
\label{eq:lc--dw-bound}
    h_\sigma\normdelim{dw}_{L^2(\tau)} \leq h_\sigma\normdelim{u}_{L^2(\tau)} + h_\sigma\normdelim{dL_\tau g}_{L^2(\tau)}.
\end{equation}
We conclude with the triangle inequality, using \eqref{eq:lc--lifting-bound} and \eqref{eq:lc--dw-bound},
\begin{align}
    \normdelim{B^k_\tau v}_{L^2(\tau)} + h_\sigma\normdelim{dB^k_\tau v}_{L^2(\tau)}
    &\leq \normdelim{w}_{L^2(\tau)} + \normdelim{(\id- \Pi)L_\tau g}_{L^2(\tau)} \\ \nonumber
    &\hspace{1em}+ h_\sigma\normdelim{dw}_{L^2(\tau)} + h_\sigma\normdelim{dL_\tau g}_{L^2(\tau)} \\ \nonumber
    &\leq Ch_\sigma\normdelim{dw}_{L^2(\tau)} + \normdelim{L_\tau g}_{L^2(\tau)} + h_\sigma\normdelim{dL_\tau g}_{L^2(\tau)} \\ \nonumber
    &\leq Ch_\sigma\normdelim{v}_{L^2(\tau)}.
\end{align}
The global estimate follows from \eqref{eq:app--claim} by summing over the top-dimensional cells.
\end{proof}

\end{proposition}

\begin{proposition}[Stability of the decomposition operator]
\label{prop:lc--stability}
Let $\mathscr H_{r,\mathrm{dec}}^\bullet$ be the decomposed homotopy operator as constructed in \cref{prop:lc--abstract-construction}, with $\mathscr H^\bullet_W$ as the $r=1$ instance of \eqref{eq:prelim--hodge-laplace}, and $B^k$ given by \eqref{eq:lc--global-B}. There is a constant $C$ depending only on $r$, $n$, $\gamma$, and the Poincar\'e constant of $\Omega$, such that for every $u\in\mathcal P_r^-\Lambda^k(\mathcal T_h)$,
\begin{equation}
    \normdelim{\mathscr H^k_{r,\mathrm{dec}}u}_{L^2(\Omega)} + \normdelim{d\mathscr H^k_{r,\mathrm{dec}}u}_{L^2(\Omega)} \leq C\normdelim{u}_{L^2(\Omega)}.
\end{equation}
\end{proposition}

\begin{proof}
Condition \eqref{eq:prelim--abstract-bstar-b} states that $\mathscr H^k_Wu\perp \mathfrak Z\mathcal P_1^-\Lambda^{k-1}(\mathcal T_h)$, and so we have a Poincar\'e inequality 
\begin{equation}
    \normdelim{\mathscr H^k_W u}_{L^2(\Omega)} \leq C\normdelim{d\mathscr H^k_Wu}_{L^2(\Omega)},
\end{equation}
with $C$ independent of $h$ (see e.g.\ \cite{Arnold-Falk-Winther:10}). Condition \eqref{eq:prelim--abstract-bstar-a} gives $d\mathscr H^k_{\mathrm W}v$ as the $L^2$-orthogonal projection of $v$ onto $d\mathcal P_1^-\Lambda^{k-1}(\mathcal T_h)$, so $\|d\mathscr H^k_{\mathrm W}\|_{L^2(\Omega)^*}= 1$. With that, and \cref{prop:lc--local-stability}, we compute
\begin{align}
    \normdelim{\mathscr H^k_{r,\mathrm{dec}}u}_{L^2(\Omega)} + \normdelim{d\mathscr H^k_{r,\mathrm{dec}}u}_{L^2(\Omega)} 
    &\leq \normdelim{\mathscr H^k_WI^ku}_{L^2(\Omega)} + \normdelim{B^k(\id - I^k)u}_{L^2(\Omega} \\ \nonumber
    &\hspace{1em}+ \normdelim{d\mathscr H^k_WI^ku}_{L^2(\Omega)} + \normdelim{dB^k(\id - I^k)u}_{L^2(\Omega)} \\ \nonumber
    &\leq C\normdelim{I^ku}_{L^2(\Omega)} + C\normdelim{(\id - I^k)u}_{L^2(\Omega)}.
\end{align}
The interpolant $I^k_\sigma$ commutes with affine pullback, because the dof in \eqref{eq:lc--whitney-interpolant} do, so $\|I^ku\|_{L^2(\sigma)}\leq C\|u\|_{L^2(\sigma)}$ for every $\sigma\in\mathcal T_h^n$ by scaling to the reference simplex. Summing over the cells gives the bound on $I^k$ and, with it, $\|(\id-I^k)u\|_{L^2(\Omega)}\leq C\|u\|_{L^2(\Omega)}$. 
\end{proof}


\section{Construction by extension over the cone}
\label{section:cone-extension}

We now give the geometric construction of a Poincar\'e operator, which is non-local, and generalises the cone construction for Whitney forms of \cite[Section~4.2]{guzman-potu:26} to arbitrary polynomial order. It corresponds to a simplified formula for the interpolation of the standard Poincar\'e operator (see \cref{section:conext--classical-relation}). To give the main idea, suppose $\Omega$ is star-shaped with respect to a point $p_0$. Then, the ``cone'' of $\sigma = [x_0, \dots, x_\ell]\in\mathcal T_h^\ell$ is $\mathbf{C}(\sigma) = [p_0, x_0, \dots, x_\ell]$. In \cite[Section 4.2]{guzman-potu:26}, a Poincar\'e operator over the Whitney forms, $\mathscr P^{k+1}:\mathcal P_1^-\Lambda^{k+1}(\mathcal T_h) \to \mathcal P_1^-\Lambda^k(\mathcal T_h)$, is defined through the canonical dof as
\begin{equation}
\label{eq:conext--poincare-lowest-order}
\int_\sigma\mathscr P^{k+1}u \coloneqq \int_{\mathbf C(\sigma)} u \qquad \forall \sigma\in \mathcal T_h^k.
\end{equation}
Using Stokes' formula along with the homotopy identity on chains $\partial \mathbf{C}(\sigma) + \mathbf{C} (\partial \sigma) = \sigma -\epsilon$ for appropriate $\epsilon$ shows immediately that the operator of \eqref{eq:conext--poincare-lowest-order} is indeed a Poincar\'e operator.

In general \cref{prop:prelim--dof-trimmed-polynomial,prop:prelim--dof-full-polynomial} let us define operators on polynomial spaces by prescribing moments against other polynomial forms. Indeed, this is exactly what is happening in \eqref{eq:conext--poincare-lowest-order}: we are prescribing the moments of $\mathscr P^{k+1}u$ against the local constant $0$-form. In the general case, in order to map a $(k{+}1)$-form $u$ to a $k$-form $\mathscr P^{k+1}u$, we must assign to each $\sigma\in \mathcal T_h^\ell$ the moments of $\mathscr P^{k+1}u$ against the appropriate $(\ell{-}k)$-forms. Note that we are not allowed to integrate $u$ against $(\ell{-}k)$-forms on $\sigma$, but we can integrate $u$ against $(\ell{-}k)$-forms on $\mathbf{C}(\sigma)$. With that observation, we consider the following interpretation of \eqref{eq:conext--poincare-lowest-order}:
\begin{equation}
\label{eq:conext--poincare-higher-order-heuristic}
\int_\sigma \mathscr P^{k+1}u\wedge q\coloneqq \int_{\mathbf{C}(\sigma)} u\wedge q,
\end{equation}
which should hold for all appropriate polynomial forms $q$, and all appropriate simplices $\sigma$. Carrying out the appropriate integration by parts produces the homotopy identity, as we verify below. What \eqref{eq:conext--poincare-higher-order-heuristic} does not yet make sense of is the right hand side, since \emph{a priori} $q$ lives only on $\sigma$. Moreover, in the case $\sigma\in\mathcal T_h^n$ (which is never needed for the Whitney forms, but is needed for the higher order spaces), $\mathbf{C}(\sigma)$ is always a degenerate $(n{+}1)$-simplex, and so this too needs to be understood correctly. We present the operator for the trimmed family $\mathcal P_r^-\Lambda^\bullet$. The construction for the degree-reducing full family is identical, using the dof of \cref{prop:prelim--dof-full-polynomial} in place of those of \cref{prop:prelim--dof-trimmed-polynomial}.

\subsection{The cone Poincar\'e operator on star-shaped domains}\label{section:conext--cone-star-shaped}

We begin by presenting the construction on domains which are star-shaped with respect to a point $p_0\in\Omega$. That is, for all $x\in\Omega$ the line segment between $p_0$ and $x$, $[p_0, x]$, lies in $\Omega$. We further suppose $p_0$ is in \emph{general position} with respect to the mesh $\mathcal{T}_h$. This means that $p_0$ is affinely independent to all $\sigma\in \mathcal T_h^\ell$, for all $\ell \in \{1, \dots, n{-}1\}$, and each cone $\mathbf C(\sigma)=[p_0, x_0, \dots, x_\ell]$ has positive $(\ell{+}1)$-dimensional measure. This assumption can be removed later with a slight modification (see \cref{rem:conext--removing-general-position}). In \cref{section:conext--general-setting,section:conext--general-construction} we discuss how the construction extends to more general domains. 

To define the operator, we now suppose we have a family of maps $\ext_\sigma:\mathcal P_{r+k-\ell-1}\Lambda^{\ell-k}(\sigma)\to \mathrsfs W_{1,1}^{\ell-k}(\mathbf{C}(\sigma))$, defined for all $k\in \{0,\dots,n\}$ and $\sigma\in\mathcal T_h^\ell$ for $k\leq \ell \leq n-1$, which satisfies
\begin{align}
\tr_{\sigma}\ext_\sigma q &=q, \tag{E.1}\label{eq:conext--E1}\\
\tr_{\mathbf{C}(\tau)}\ext_\sigma q&=\ext_\tau\tr_\tau q \qquad \forall\tau\triangleleft\sigma,\, k\ge 1,\tag{E.2}\label{eq:conext--E2}\\
d\,\ext_\sigma q&=\ext_\sigma dq. \tag{E.3}\label{eq:conext--E3}
\end{align}
We construct such extensions explicitly in \cref{section:conext--extensions}.

To deal with the dofs associated to top-dimensional cells, we use locally defined operators $K_\sigma: \mathcal P_{r+k-n-1}\Lambda^{n-k}(\sigma) \to \mathcal P_{r+k-n}\Lambda^{n-k-1}(\sigma)$, one for each $\sigma\in\mathcal T_h^n$, such that $dK_\sigma + K_\sigma d = \id-\pi_0$ where $\pi_0 q = q(a)$ if $q$ is a 0-form and $\pi_0 = 0$ otherwise. In particular, the scaled Koszul operator suffices (see e.g. \cite{Arnold-Falk-Winther:06}). Given $q\in \mathcal P_{r+k-n-1}\Lambda^{n-k}(\sigma)$, write $q=\sum_j q_j$, where $q_j$ has homogeneous coefficients of degree $j$, with respect to $x-a$ for some point $a$ in $\sigma$. We then take
\begin{equation}
K_{\sigma} q \coloneqq \sum_{j} \frac{\kappa_{a} q_j}{j+n-k}.
\end{equation}
Note that this is a purely algebraic operator.

\begin{definition}[Poincar\'e operator from the cone]
\label{def:conext--poincare}
The discrete Poincar\'e operator is defined by setting the dofs as follows. For all $\sigma\in  \mathcal T_h^\ell$ with $\ell\leq n{-}1$, and $q\in \mathcal P_{r+k-\ell-1}\Lambda^{\ell-k}(\sigma)$, we set
\begin{subequations}
\label{eq:conext--poincare-def-subeq}
\begin{equation}
\label{eq:conext--poincare-def}
\int_\sigma \mathscr P_{r, \mathrm{co}}^{k+1} u\wedge q \coloneqq \int_{\mathbf{C}(\sigma)} u\wedge \ext_\sigma q.
\end{equation}
For $\sigma\in\mathcal T_h^n$ we instead set, for all $q\in \mathcal P_{r+k-n-1}\Lambda^{n-k}(\sigma)$,
\begin{equation}
\label{eq:conext--top-cell}
    \int_\sigma\mathscr P_{r,\mathrm{co}}^{k+1}u\wedge q 
    \coloneqq (-1)^{k+1}\int_{\sigma} u \wedge K_{\sigma} q 
    +(-1)^{k}\sum_{\tau \triangleleft \sigma} o_{\sigma\tau} \int_{\mathbf C(\tau)} u \wedge \ext_\tau (\tr_\tau K_{\sigma} q).
\end{equation}
\end{subequations}
\end{definition}

By \cref{prop:prelim--dof-trimmed-polynomial} the dof appearing in \eqref{eq:conext--poincare-def}-\eqref{eq:conext--top-cell} are unisolvent on $\mathcal P_r^-\Lambda^k(\mathcal T_h)$, so \eqref{eq:conext--poincare-def}--\eqref{eq:conext--top-cell} determines a unique element of that space. For the Whitney forms $\mathcal P_1^-\Lambda^k$ the only test forms are the constants on $k$-simplices, and \cref{def:conext--poincare} reduces precisely to the lowest order operator \eqref{eq:conext--poincare-lowest-order}.

\begin{theorem}[Homotopy formula]
\label{thm:conext--homotopy-formula}
For $1\leq k\leq n$ and $u\in \mathcal P_r^-\Lambda^k(\mathcal T_h)$, it holds
\begin{equation}
\label{eq:conext--homotopy-formula}
d\mathscr P_{r, \mathrm{co}}^k u+\mathscr P_{r, \mathrm{co}}^{k+1} du = u,
\end{equation}
whilst for $u\in\mathcal P_r^-\Lambda^0(\mathcal T_h)$,
\begin{equation}
\label{eq:conext--homotopy-formula-0}
\mathscr P_{r, \mathrm{co}}^1du=u - u(p_0).
\end{equation}  
\end{theorem}

To prove \cref{thm:conext--homotopy-formula} we use the following pair of lemmas.
\begin{lemma}[Stokes' theorem on the cone]\label{lem:conext--cone-stokes}
Assume \eqref{eq:conext--E1}--\eqref{eq:conext--E3} hold. Let $\sigma\in\mathcal T_h^\ell$ with $\ell\le n{-}1$, $q\in \mathcal P_{r+k-\ell-1}\Lambda^{\ell-k}(\sigma)$, and $u\in \mathrsfs W^k_{\infty,\infty}(\Omega)$. Then $u\wedge\ext_\sigma q$ and $d(u\wedge\ext_\sigma q)$ are integrable on $\mathbf C(\sigma)$, and when $k\ge 1$,
\begin{subequations}
\label{eq:conext--stokes-cone}
    \begin{equation}
    \label{eq:conext--stokes-cone-a}
    \int_{\mathbf C(\sigma)} d(u\wedge\ext_\sigma q)
    =\int_\sigma u\wedge q-\sum_{\tau\triangleleft \sigma} o_{\sigma\tau}\int_{\mathbf C(\tau)}u\wedge\ext_\tau\tr_\tau q.
    \end{equation}
For $k=0$ we instead get
    \begin{equation}
    \label{eq:conext--stokes-cone-b}
    \int_{\mathbf C(\sigma)} du\wedge\ext_\sigma q
    =\int_\sigma u\wedge q
    -u(p_0)\int_\sigma q.
    \end{equation}
\end{subequations}
\end{lemma}

\begin{proof}
As $\ext_\sigma q$ is of regularity $\mathrsfs W_{1,1}^{\ell-k}$, and $u$ is of regularity $\mathrsfs W_{\infty,\infty}^k$ on a bounded domain, it holds that $u\wedge \ext_\sigma q\in\mathrsfs W_{1,1}^\ell(\mathbf C(\sigma))$. Therefore we can apply Stokes' theorem, giving us
\begin{equation}
\label{eq:conext--cone-stokes-direct}
    \int_{\mathbf C(\sigma)} d(u\wedge \ext_\sigma q) 
    = \int_\sigma u\wedge q - \sum_{\tau\triangleleft \sigma} o_{\sigma\tau} \int_{\mathbf C(\tau)} u\wedge \tr_{\mathbf C(\tau)} \ext_\sigma q.
\end{equation}
Here we used \eqref{eq:conext--E1} and the fact that $\partial \mathbf C(\sigma) = \sigma-\sum_\tau o_{\sigma\tau}\mathbf C(\tau)$. For $k\geq 1$ the identity \eqref{eq:conext--stokes-cone-a} now follows immediately from \eqref{eq:conext--E3}. The remaining case is shown by an induction argument on $\ell=\dim \sigma$.

Let $k=0$. For $\ell=0$ we have that $\sigma =[x_0]$ is a point, and $\mathbf C(\sigma)$ is the straight line $[p_0,x_0]$. Then $q$ is a constant, and $\ext_\sigma q\equiv q$ by \eqref{eq:conext--E1} and \eqref{eq:conext--E3}. Thus \eqref{eq:conext--stokes-cone-b} states $\int_{[p_0,x_0]}q\,du=q\,(u(x_0)-u(p_0))$, which is the fundamental theorem of calculus. The point evaluation is well defined since $u$ is Lipschitz. 

Let now $\ell\ge1$, and suppose \eqref{eq:conext--stokes-cone-b} holds for all facets $\tau\triangleleft \sigma$. As $q\in\mathcal P_{r-\ell-1}\Lambda^\ell(\sigma)$ is of top form degree on $\sigma$, there exists $\omega_q\in\mathcal P_{r-\ell}\Lambda^{\ell-1}(\sigma)$ such that $q=d\omega_q$. Applying \eqref{eq:conext--stokes-cone-b} to the pair of $du$ and $\omega_q$, and using \eqref{eq:conext--E3}, we find that
\begin{align}
\label{eq:conext--cone-stokes--proof-b-1}
    \int_{\mathbf C(\sigma)} du \wedge \ext_\sigma q
    &= -\int_{\mathbf C(\sigma)} d(du\wedge \ext_\sigma \omega_q) \\ \nonumber
    &= -\int_\sigma du \wedge \omega_q + \sum_{\tau\triangleleft \sigma} o_{\sigma\tau} \int_{\mathbf C(\tau)} du \wedge \ext_\tau\tr_\tau \omega_q
\end{align}
We may then invoke the induction hypothesis to rewrite the right-most summands, which looks like
\begin{equation}
\label{eq:conext--cone-stokes--proof-b-2}
    \int_{\mathbf C(\tau)} du \wedge \ext_\tau\tr_\tau \omega_q
    = \int_\tau u \wedge \tr_\tau \omega_q - u(p_0)\int_\tau \tr_\tau \omega_q.
\end{equation}
Taking \eqref{eq:conext--cone-stokes--proof-b-2} to the sum over all facets, Stokes' theorem yields
\begin{equation}
\label{eq:conext--cone-stokes--proof-b-3}
    \sum_{\tau\triangleleft \sigma}o_{\sigma\tau}\int_{\mathbf C(\tau)} du \wedge \ext_\tau\tr_\tau \omega_q
    = \int_\sigma du\wedge \omega_q + \int_\sigma u\wedge q - u(p_0)\int_\sigma q.
\end{equation}
Finally, the desired equality \eqref{eq:conext--stokes-cone-b} follows by plugging \eqref{eq:conext--cone-stokes--proof-b-3} into \eqref{eq:conext--cone-stokes--proof-b-1}. This completes the induction and the proof.
\end{proof}

\begin{lemma}[Integration by parts on the degenerate cone]
\label{lem:conext--degenerate-stokes}

Assume \eqref{eq:conext--E1}--\eqref{eq:conext--E3} hold. Let $\sigma\in\mathcal T_h^n$, $q\in \mathcal P_{r+k-n-1}\Lambda^{n-k-1}(\sigma)$, and $u\in \mathrsfs W^k_{\infty,\infty}(\Omega)$. Then the below integrals are well-defined, and for $k\geq  1$ it holds
\begin{subequations}
\label{eq:conext--degenerate-stokes}
\begin{equation}
\label{eq:conext--degenerate-stokes-a}
    \int_\sigma du \wedge q - \sum_{\tau\triangleleft \sigma} o_{\sigma\tau} \int_{\mathbf C(\tau)} du\wedge \ext_\tau\tr_\tau q
    = (-1)^{k+1}\int_\sigma u\wedge dq + (-1)^k \sum_{\tau\triangleleft \sigma} o_{\sigma\tau} \int_{\mathbf C(\tau)} u\wedge \ext_\tau\tr_\tau dq.
\end{equation}
For $k=0$ we instead get
\begin{equation}
\label{eq:conext--degenerate-stokes-b}
    \int_\sigma du \wedge q - \sum_{\tau\triangleleft \sigma} o_{\sigma\tau} \int_{\mathbf C(\tau)} du\wedge \ext_\tau\tr_\tau q
    = -\int_\sigma u\wedge dq + u(p_0) \int_\sigma dq.
\end{equation}
\end{subequations}
\end{lemma}

\begin{proof}
    The strategy is the same as in the proof of \cref{lem:conext--cone-stokes}, except that we avoid integrating over $\mathbf C(\sigma)$, since it is degenerate. The key to the argument lies in the identity
    \begin{equation}
    \label{eq:conext--degenerate-stokes--proof-1}
        \sum_{\tau\triangleleft \sigma}\sum_{\tau'\triangleleft \tau} o_{\sigma\tau}o_{\tau\tau'} = 0,
    \end{equation}
    which follows from the fact that $\partial\circ \partial=0$. Applying \cref{lem:conext--cone-stokes} on $\mathbf C(\tau)$ for facets $\tau \triangleleft \sigma$ then yields, for $k\geq 1$,
    \begin{align}
    \label{eq:conext--degenerate-stokes--proof-2}
        \sum_{\tau\triangleleft \sigma} \int_{\mathbf C(\tau)} du \wedge \ext_\tau\tr_\tau q 
        &= \sum_{\tau\triangleleft \sigma} \int_{\mathbf C(\tau)} d(u \wedge \ext_\tau\tr_\tau q) 
        + (-1)^{k+1}\sum_{\tau\triangleleft \sigma} \int_{\mathbf C(\tau)} u \wedge \ext_\tau\tr_\tau dq \\ \nonumber
        &= \int_{\partial\sigma} u\wedge q 
        + (-1)^{k+1}\sum_{\tau\triangleleft \sigma} \int_{\mathbf C(\tau)} u \wedge \ext_\tau\tr_\tau dq.
    \end{align}
    Indeed, the weighted sum of the boundary terms that appear in \eqref{eq:conext--stokes-cone-a} vanishes by \eqref{eq:conext--degenerate-stokes--proof-1} and \eqref{eq:conext--E2} (note that $q$ is never of top form degree, so we always have \eqref{eq:conext--E2}). For $k=0$, \cref{lem:conext--cone-stokes} instead gives
    \begin{equation}
    \label{eq:conext--degenerate-stokes--proof-3}
        \sum_{\tau\triangleleft \sigma} \int_{\mathbf C(\tau)} du \wedge \ext_\tau\tr_\tau q 
        = \int_{\partial\sigma} u\wedge q 
        - u(p_0) \int_{\partial\sigma} q.
    \end{equation}
    Finally, a simple application of Stokes' theorem yields
    \begin{equation}
    \label{eq:conext--degenerate-stokes--proof-4}
        \int_\sigma du\wedge q = \int_{\partial\sigma} u\wedge q + (-1)^{k+1}\int_\sigma u\wedge dq,
    \end{equation}
    and the statement \eqref{eq:conext--degenerate-stokes} follows by subtracting either \eqref{eq:conext--degenerate-stokes--proof-2} or \eqref{eq:conext--degenerate-stokes--proof-3}, depending on $k$, from \eqref{eq:conext--degenerate-stokes--proof-4}.
\end{proof}

\begin{proof}[Proof of \cref{thm:conext--homotopy-formula}]
By \cref{prop:prelim--dof-trimmed-polynomial} it suffices to test both sides against $q\in\mathcal P_{r+k-\ell-1}\Lambda^{\ell-k}(\sigma)$ for every $\sigma\in\mathcal T_h^\ell$ and every $\ell\in\set{k,\dots,n}$. We split it up into four cases, depending on whether or not $k=0$, and whether $\ell=n$. Suppose first that $k\geq1$ and $\ell \leq n{-}1$. Since $\mathscr P^k_{r,\mathrm{co}}u\wedge q$ is an $(\ell{-}1)$-form on the $\ell$-simplex $\sigma$, Stokes' theorem on $\sigma$ gives
\begin{equation}
\label{eq:conext--homotopy-proof-term10}
\int_\sigma d\mathscr P_{r, \mathrm{co}}^k u \wedge q
= \int_{\partial \sigma}\mathscr P_{r, \mathrm{co}}^ku\wedge q + (-1)^k\int_\sigma\mathscr P_{r, \mathrm{co}}^k u\wedge dq.
\end{equation}
The two terms on the right are again of the form \eqref{eq:conext--poincare-def}, one degree lower. Indeed, for a facet $\tau\triangleleft\sigma$ we have $\tr_\tau q\in \mathcal P_{r+(k-1)-\dim\tau-1}\Lambda^{\dim\tau-(k-1)}(\tau)$ by \cref{prop:prelim--trace-polynomial}, and $dq\in\mathcal P_{r+(k-1)-\ell-1}\Lambda^{\ell-(k-1)}(\sigma)$ by \cref{prop:prelim--extder-polynomial}, meaning both are admissible test forms for $\mathscr P^k_{r,\mathrm{co}}$. Hence \eqref{eq:conext--homotopy-proof-term10} becomes
\begin{equation}
\label{eq:conext--homotopy-proof-term11}
\int_\sigma d\mathscr P_{r, \mathrm{co}}^k u \wedge q
= \sum_{\tau\triangleleft \sigma} o_{\sigma\tau} \int_{\mathbf C(\tau)} u\wedge \ext_{\tau}\tr_\tau q
+ (-1)^k\int_{\mathbf C(\sigma)} u\wedge \ext_\sigma dq .
\end{equation}
Moreover $du\in\mathcal P_r^-\Lambda^{k+1}(\mathcal T_h)$ by \cref{prop:prelim--extder-polynomial}, so \eqref{eq:conext--poincare-def} gives
\begin{equation}
\label{eq:conext--homotopy-proof-term21}
\int_\sigma \mathscr P_{r, \mathrm{co}}^{k+1} du\wedge q = \int_{\mathbf C(\sigma)}du\wedge \ext_\sigma q .
\end{equation}
Adding \eqref{eq:conext--homotopy-proof-term11} to \eqref{eq:conext--homotopy-proof-term21} and applying \cref{lem:conext--cone-stokes} yields $\int_\sigma\bigl(d\mathscr P^k_{r,\mathrm{co}}u+\mathscr P^{k+1}_{r,\mathrm{co}}du\bigr)\wedge q=\int_\sigma u\wedge q$, as desired. If instead $k=0$, then $\int_\sigma \mathscr P_{r,\mathrm{co}}^1du \wedge q=\int_\sigma (u-u(p_0))\wedge q$ follows immediatly from \eqref{eq:conext--stokes-cone-b}.

It remains to check the cases where $\ell=n$, so suppose $\sigma\in\mathcal T_h^n$, and recall that the dofs on top-dimensional cells are given by \eqref{eq:conext--top-cell}. Consider first the case $k\geq 1$. Plugging the new definition \eqref{eq:conext--top-cell} into the right-most term of \eqref{eq:conext--homotopy-proof-term10}, which remains valid, we find that
\begin{equation}
\label{eq:conext--homotopy-proof-term31}
    \int_\sigma d\mathscr P_{r,\mathrm{co}}^ku\wedge q
    = \sum_{\tau\triangleleft \sigma} o_{\sigma\tau} \int_{\mathbf C(\tau)} u \wedge \ext_\tau \tr_\tau q
    + \int_\sigma u\wedge K_\sigma dq
    - \sum_{\tau\triangleleft \sigma} o_{\sigma\tau}\int_{\mathbf C(\tau)} u\wedge \ext_\tau(\tr_\tau Kdq).
\end{equation}
Note that the boundary terms in \eqref{eq:conext--homotopy-proof-term10} still adhere to the same definition as before, \eqref{eq:conext--poincare-def}. The other term is computed using \cref{lem:conext--degenerate-stokes} as follows,
\begin{align}
\label{eq:conext--homotopy-proof-term32}
    \int_{\sigma}\mathscr P_{r,\mathrm{co}}^{k+1}du\wedge q
    &= (-1)^{k+1}\int_\sigma du \wedge K_\sigma q + (-1)^k \sum_{\tau\triangleleft \sigma} o_{\sigma\tau} \int_{\mathbf C(\tau)} du \wedge \ext_\tau (\tr_\tau K_\sigma q) \\ \nonumber
    &= \int_\sigma u \wedge dK_\sigma q - \sum_{\tau\triangleleft\sigma} o_{\sigma\tau} \int_{\mathbf C(\tau)} u\wedge \ext_\tau(\tr_\tau dK_\sigma q).
\end{align}
Since $q$ is never a $0$-form, it always holds that $dK_\sigma q+K_\sigma dq=q$. Thus upon adding \eqref{eq:conext--homotopy-proof-term31} to \eqref{eq:conext--homotopy-proof-term32}, the boundary terms cancel, and we are left with $\int_\sigma\bigl(d\mathscr P^k_{r,\mathrm{co}}u+\mathscr P^{k+1}_{r,\mathrm{co}}du\bigr)\wedge q=\int_\sigma u\wedge q$. Finally, for $k=0$ we have that $q$ is of top form degree, and so $q=dK_\sigma q$. Taken together with \cref{lem:conext--degenerate-stokes} this gives us
\begin{align}
    \int_\sigma \mathscr P_{r,\mathrm{co}}^1du \wedge q 
    &= -\int_\sigma du\wedge K_\sigma q 
    + \sum_{\tau\triangleleft \sigma} o_{\sigma\tau} \int_{\mathbf C(\tau)} du \wedge \ext_\tau (\tr_\tau K_\sigma q) \\ \nonumber
    &= \int_\sigma u\wedge q
    - u(p_0)\int_\sigma q,
\end{align}
as desired.
\end{proof}

\begin{corollary}
\label{cor:conext--is-poincare}
$\mathscr P^\bullet_{r,\mathrm{co}}$ is a Poincar\'e operator on $\mathcal P_r^-\Lambda^\bullet(\mathcal T_h)$ in the sense of \cref{def:prelim--poincare}, with cohomology projection $\pi^0u=u(p_0)$ and $\pi^k=0$ for $k>0$.
\end{corollary}

\begin{remark}[Removing the general position assumption on $p_0$]\label{rem:conext--removing-general-position}

Suppose $\sigma\in\mathcal T_h^\ell$ is \emph{cone degenerate} in the sense that $p_0$ is affinely dependent on $\sigma$. Then the dof associated to $\sigma$ are set to zero and the homotopy identity \eqref{eq:conext--homotopy-formula} no longer holds on $\sigma$. To fix this we can do as we did for the top-dimensional cells, and use the Koszul operator. We keep \eqref{eq:conext--poincare-def} for the cone non-degenerate faces $\tau\triangleleft\sigma$, and for $q\in\mathcal P_{r+k-\ell-1}\Lambda^{\ell-k}(\sigma)$ we set
\begin{equation}\label{eq:conext--degenerate-cell}
\int_\sigma\mathscr P^{k+1}_{r,\mathrm{co}}u\wedge q
\coloneqq(-1)^{k+1}\int_\sigma u\wedge K_{p_0}q
+(-1)^k\sum_{\tau\triangleleft\sigma}o_{\sigma\tau}\int_{\mathbf C(\tau)}u\wedge\ext_\tau(\tr_\tau K_{p_0}q).
\end{equation}
With this modification, \Cref{thm:conext--homotopy-formula} remains true. For cone non-degenerate simplices, the proof is unchanged, because \cref{lem:conext--cone-stokes} only involves the cones over $\sigma$ and its facets. Otherwise, the verification is essentially the argument given for top-dimensional cells, with $n$ replaced by $\dim\sigma$ and by using induction over $\dim\sigma$ (the base case is immediate since $\mathbf{C}(x)$ is not degenerate for all $x\neq p_0$). In particular, the proof uses only that $dK_{p_0}+K_{p_0}d=\id-\pi_0$ and the homotopy identity inductively established on the facets.
\end{remark}

\subsection{The extension operators on star-shaped domains}\label{section:conext--extensions}
Much work has been done on developing commuting extension operators for finite elements \cite{DGS:08, DGS:09, DGS:12}. These extensions were developed for extending polynomial traces from the boundary of a tetrahedron into its interior for the analysis of high order finite elements. However they do not in general satisfy \eqref{eq:conext--E1}--\eqref{eq:conext--E3}. In particular, the conditions \eqref{eq:conext--E1}-\eqref{eq:conext--E3} necessarily introduce some singularities at $p_0$ which immediately disallow the use of polynomial extensions.

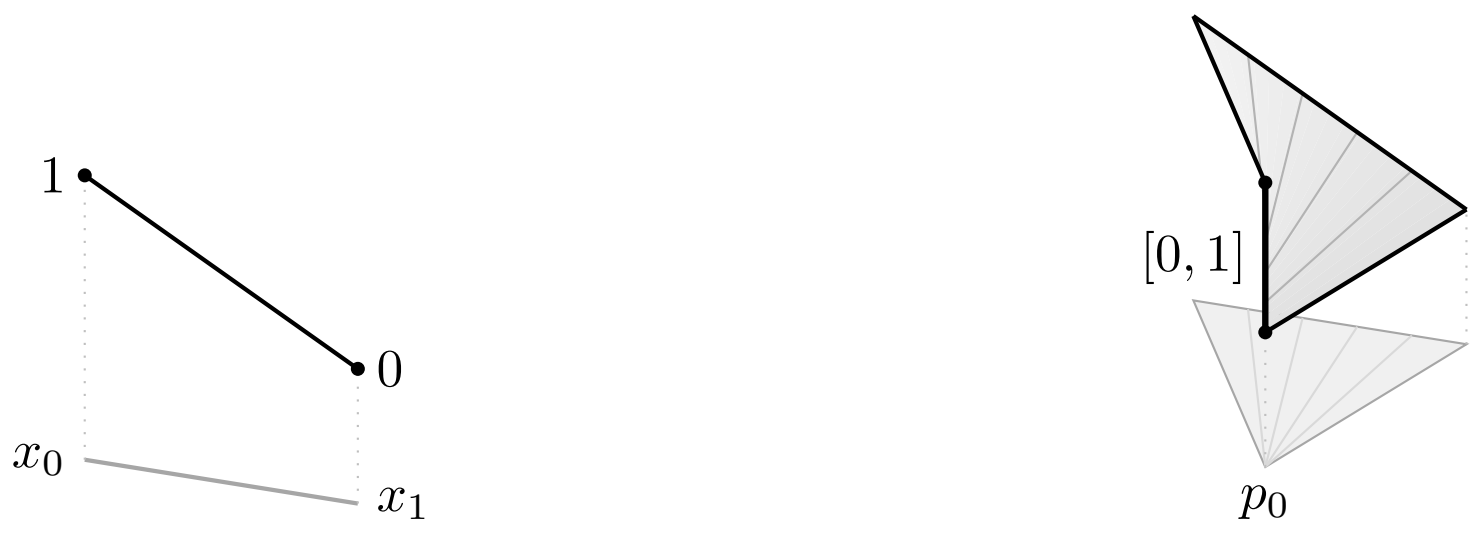
\begin{figure}[ht]
\centering
\begin{subfigure}[b]{0.30\linewidth}\centering
\begin{tikzpicture}[x={(0.966cm,-0.155cm)},y={(0.259cm,0.580cm)},z={(0cm,0.90cm)},scale=1.10]
\draw[gray!70,thick] (0,1.5,-0.9)--(1.7,1.5,-0.9);
\draw[gray!50,dotted] (0,1.5,-0.9)--(0,1.5,1);
\draw[gray!50,dotted] (1.7,1.5,-0.9)--(1.7,1.5,0);
\draw[thick] (0,1.5,1)--(1.7,1.5,0);
\fill (0,1.5,1) circle (1.2pt); \fill (1.7,1.5,0) circle (1.2pt);
\node[left] at (0,1.5,1) {$1$};
\node[right] at (1.7,1.5,0) {$0$};
\node[left] at (0,1.5,-0.9) {$x_0$};
\node[right] at (1.7,1.5,-0.9) {$x_1$};
\end{tikzpicture}
\caption{$q$ on $\sigma$}\label{fig:conext--graph-q}
\end{subfigure}\hspace{2cm}
\begin{subfigure}[b]{0.34\linewidth}\centering
\begin{tikzpicture}[x={(0.966cm,-0.155cm)},y={(0.259cm,0.580cm)},z={(0cm,0.90cm)},scale=1.10]
\fill[gray!12] (0.85,0,-0.9)--(0,1.5,-0.9)--(1.7,1.5,-0.9)--cycle;
\draw[gray!70] (0.85,0,-0.9)--(0,1.5,-0.9)--(1.7,1.5,-0.9)--cycle;
\foreach \s in {0.2,0.4,0.6,0.8} {\draw[gray!30] (0.85,0,-0.9)--({1.7*\s},1.5,-0.9);}
\foreach \i in {0,...,17}{
  \pgfmathsetmacro{\sa}{\i/18}\pgfmathsetmacro{\sb}{(\i+1)/18}
  \pgfmathsetmacro{\gr}{10+14*\sa}
  \fill[gray!\gr] (0.85,0,{1-\sa})--({1.7*\sa},1.5,{1-\sa})--({1.7*\sb},1.5,{1-\sb})--(0.85,0,{1-\sb})--cycle;}
\foreach \s in {0.2,0.4,0.6,0.8} {\draw[gray!60] (0.85,0,{1-\s})--({1.7*\s},1.5,{1-\s});}
\draw[gray!50,dotted] (0.85,0,-0.9)--(0.85,0,0);
\draw[gray!50,dotted] (1.7,1.5,-0.9)--(1.7,1.5,0);
\draw[thick] (0.85,0,1)--(0,1.5,1);
\draw[thick] (0.85,0,0)--(1.7,1.5,0);
\draw[thick] (0,1.5,1)--(1.7,1.5,0);
\draw[very thick] (0.85,0,0)--(0.85,0,1);
\fill (0.85,0,1) circle (1.2pt); \fill (0.85,0,0) circle (1.2pt);
\node[left] at (0.85,0,0.5) {$[0,1]$};
\node[below] at (0.85,0,-0.9) {$p_0$};
\end{tikzpicture}
\caption{$\ext_\sigma q$ on $\mathbf{C}(\sigma)$}\label{fig:conext--graph-cone}
\end{subfigure}
\caption{Graph of the test form $q=\lambda_0$ and of its extension to the geometric cone, for $\sigma=[x_0,x_1]$. The induced extension on the geometric cone has a singularity at $p_0$ as depicted in (B).}
\label{fig:conext--graphs}
\end{figure}

Indeed, let $\sigma=[x_0,x_1]$ be an edge, and suppose the cone $[p_0,x_0,x_1]$ is non-degenerate. We consider the extension of the barycentric coordinate $\lambda_0$ associated to $x_0$. Using \eqref{eq:conext--E1} and \eqref{eq:conext--E2} yields
\begin{equation}
    \tr_{\mathbf{C}(x_i)} \ext_\sigma \lambda_0 = \ext_{x_i} \tr_{x_i} \lambda_0 = \begin{cases}
        \ext_{x_0}1 & \text{if $i=0$},\\
        \ext_{x_1}0 & \text{if $i=1$}.
    \end{cases}
\end{equation}
Then \eqref{eq:conext--E3} together with \eqref{eq:conext--E1} says that the extension of a constant is itself, so that the right-hand side is $\delta_{i0}$. Thus the value of $\ext_\sigma^0\lambda_0$ is ambiguous on the intersection of $\mathbf{C}(x_0)$ with $\mathbf{C}(x_1)$, i.e. there is a singularity at this corner (see \cref{fig:conext--graphs}). By using rational forms we are able to construct explicit extensions satisfying the necessary assumptions \eqref{eq:conext--E1}-\eqref{eq:conext--E3}.

\begin{definition}[Rational extension operators]\label{def:conext--rational-extensions}
Let $\lambda_0,\dots,\lambda_\ell$ be the barycentric coordinates of $\sigma$ and $\mu_{p_0},\mu_0,\dots,\mu_{\ell}$
those of $[p_0,x_0,\dots,x_\ell]$, with $\mu_i$ associated to $x_i$. Let $\rho_\sigma:\mathbf C(\sigma)\setminus\{p_0\}\to\sigma$ be the radial retraction, $\rho_\sigma\bigl((1-t)p_0+ty\bigr)\coloneqq y$ for $y\in\sigma$ and $t\in(0,1]$. Then the rational extension operator is defined by $\ext_\sigma q\coloneqq \rho_\sigma^* q$. 
Upon writing $x=(1-t)p_0+ty$ we find that
\begin{equation}
    \ext_\sigma^0 \lambda_i(x) = \lambda_i(y) = \frac{1}{t}\mu_i(x) - \frac{1-t}{t}\underbrace{\mu_i(p_0)}_{=\,0}.
\end{equation}
One can then check that $t=1-\mu_{p_0}(x)$, so we end up with
\begin{equation}
    \ext_\sigma^0\lambda_i = \frac{\mu_i}{1-\mu_{p_0}} \eqqcolon \nu_i,\qquad i=0,\dots,\ell.
\end{equation}
Since pullbacks distribute over the wedge and commute with $d$, it holds
\begin{equation}
\label{eq:conext--radial-ext-explicit}
\ext_\sigma^k\bigl(\lambda^\alpha\,d\lambda_{i_1}\wedge\cdots\wedge d\lambda_{i_k}\bigr)
=\nu^\alpha\,d\nu_{i_1}\wedge\cdots\wedge d\nu_{i_k},
\qquad
d\nu_i=\frac{(1-\mu_{p_0})\,d\mu_i+\mu_i\,d\mu_{p_0}}{(1-\mu_{p_0})^2}.
\end{equation}
\end{definition}
An illustration of this extension is given in Figure~\ref{fig:conext--graphs}.

\begin{proposition}[Properties of the rational extensions]\label{lem:conext--rational-extensions}
Let $\sigma \in \mathcal T_h^\ell$ and $q\in \mathcal P_m\Lambda^k(\sigma)$ for any $m\ge 0$. Then $\ext_\sigma^k q \in \mathrsfs W_{1,1}^k(\mathbf{C}(\sigma))$ and satisfies \eqref{eq:conext--E1}-\eqref{eq:conext--E3}.
\end{proposition}

\begin{proof}
The properties \eqref{eq:conext--E1}-\eqref{eq:conext--E3} are immediate from the fact that $\ext_\sigma q\coloneqq \rho_\sigma^* q$, since $d$ commutes with pullbacks, and the trace maps are themselves pullbacks of the appropriate inclusions. One can also verify these properties directly using the formula \eqref{eq:conext--radial-ext-explicit}.

It remains to show that $\ext_\sigma^k q \in \mathrsfs W_{1,1}^k(\mathbf{C}(\sigma))$. To see this, we use the formula \eqref{eq:conext--radial-ext-explicit}. First, we let $t = 1-\mu_{p_0} = \mu_1 + \cdots + \mu_k$. Then, $\mu_i\le t$, and so we have that $0\le\nu_i\le1$ on $\mathbf C(\sigma)$. Also, 
\begin{equation}
|d\nu_i|=\frac{|t\,d\mu_i+\mu_i\,d\mu_{p_0}|}{t^2}\le\frac{|d\mu_i|+|d\mu_{p_0}|}{t},
\end{equation}
where $|d\mu_i|$ and $|d\mu_{p_0}|$ are constants. From these two pointwise bounds, we see that $|\ext_\sigma q|\le C_q\,t^{-k}$ for $q\in \mathcal P_m\Lambda^k(\sigma)$. Moreover, since $t$ is affine with $t(p_0)=0$ and $t=1$ on $\sigma$, it holds that the sublevel set $\theta_s \coloneqq \{x\in\mathbf C(\sigma):t(x)\le s\}$ has measure $\abs{\theta_s} = s^{\ell+1}|\mathbf C(\sigma)|$ for $0<s\le1$. Thus, we compute
\begin{equation}
\int_{\mathbf C(\sigma)}t^{-k}\,dx
=\int_0^1 s^{-k}\,\frac{d}{ds}\abs{\theta_s}\,ds
= (\ell+1)\,|\mathbf C(\sigma)|\int_0^1 s^{\ell-k}\,ds,
\end{equation}
which is finite precisely when $k\leq \ell$. Consequently $\ext_\sigma^k q\in L^1\Lambda^k(\mathbf C(\sigma))$. By \eqref{eq:conext--E3}, $d\,\ext_\sigma q=\ext_\sigma dq$, which lies in $L^1\Lambda^{k+1}(\mathbf C(\sigma))$ by the same argument if $k<\ell$, and vanishes if $k=\ell$. We then have that $\ext_\sigma^k q\in\mathrsfs W^k_{1,1}(\mathbf C(\sigma))$, which completes the proof.
\end{proof}

As noted in \cite{guzman-potu:26}, the right hand side of \eqref{eq:conext--poincare-def} involves integrating over all mesh cells that intersect $\mathbf{C}(\sigma)$ which may be challenging in practice. Moreover, using the above extensions, one has to integrate rational forms, which introduces another difficulty. We note that the use of rational functions in finite elements goes back to the work of Wachspress \cite{Wachspress:75} and Zienkiewicz \cite{Zienkiewicz:71}. Moreover, efficient methods for integrating rational functions have been recently studied \cite{Diening-Storn-Tscherpel:25, Storn:26}. A particularly relevant recent work involving rational functions are the ``blow-up Whitney forms'' studied in \cite{BK-Gawlik:25}, which were originally developed in \cite{BGM:1991}.

\begin{remark}[Relation to blow-up Whitney forms]
\label{rem:conext--blowup}
    The extrusion $[0,1]\times \sigma$ is a ``blow-up'' of the simplex $[p_0, x_0, \dots, x_\ell]$ at the vertex $p_0$ and $(t,\lambda_0, \dots, \lambda_\ell)$ are quasi-cylindrical coordinates in the sense of \cite[Section 2.1]{BK-Gawlik:25}. In \cite{BK-Gawlik:25, BGM:1991}, a simplex is blown up along all its lower-dimensional faces, and Whitney forms on the blow-up are pushed forward to rational forms on the simplex. Thus, the lowest order instances of the extensions used here are, in a sense, a special case of the blow-up construction of \cite{BK-Gawlik:25} where a single vertex is blown up and the form is pulled back from the base of the resulting prism. For instance, with $\sigma = [x_0, x_1]$ and $q=\lambda_0$, \eqref{eq:conext--radial-ext-explicit} gives $\ext_\sigma q = \mu_0/(\mu_0+\mu_1)$ which are exactly the same shape as the blow-up Whitney forms of \cite{BK-Gawlik:25}.
\end{remark}

\subsection{Example of the construction for star-shaped \texorpdfstring{$\Omega\subset\mathbb R^3$}{3d domains}}
\label{section:cone3d}

We now write out \cref{def:conext--poincare} in vector proxies for the operator that produces curl potentials in three dimensions
\begin{equation}
\label{eq:conext--homotopy-curl-div}
\div\mathbf u=0\quad\implies\quad\mathbf u=\curl\mathscr P^2_{r,\mathrm{co}}\mathbf u.
\end{equation}
We denote by $\mathsf{RT}_r(\mathcal T_h)$ the Raviart--Thomas space of order $r$ and by $\mathsf{N1}_r(\mathcal T_h)$ the first-kind N\'ed\'elec space of order $r$, which are the proxies of $\mathcal P_r^-\Lambda^2(\mathcal T_h)$ and $\mathcal P_r^-\Lambda^1(\mathcal T_h)$ respectively.

As before, for $\sigma\in\mathcal T_h$, either an edge or a face, let $\mu_{p_0},\mu_0,\cdots,\mu_\ell$ denote the barycentric coordinates of $\mathbf C(\sigma)$. Recall that for a multi-index $\alpha$ we write $\lambda^\alpha\coloneqq\lambda_0^{\alpha_0}\cdots\lambda_\ell^{\alpha_\ell}$ and $\abs{\alpha}\coloneqq\sum_{i=0}^\ell\alpha_i$. Then for a scalar function $q$ on an edge $e$, respectively for a tangential field $\mathbf q$ on a face $f$, we have the following expressions for the rational extension
\begin{subequations}
    \begin{alignat}{2}\label{eq:cone3d--ext}
        \ext_e^0 \,\lambda^\alpha &= \frac{\mu^\alpha}{(1-\mu_{p_0})^{\abs{\alpha}}},\\
        \ext_f^1 (\lambda^\alpha\nabla \lambda_i) &= \frac{ \mu^\alpha\nabla \mu_i}{(1-\mu_{p_0})^{\abs{\alpha}+1}} + \frac{\mu^\alpha \mu_i \nabla \mu_{p_0}}{(1-\mu_{p_0})^{\abs{\alpha}+2}}.
    \end{alignat}
\end{subequations}

We orient each edge $e$ by a unit tangent $\mathbf t_e$, and denote by $\boldsymbol \nu_{\mathbf C(e)}$ the unit normal of the geometric cone $\mathbf C(e)$ pointing in the direction of $(y-p_0)\times\mathbf t_e$, where $y\in e$. Each face $f$ is oriented by its unit normal $\boldsymbol\nu_f$ pointing away from $p_0$, and $\boldsymbol{\nu}_T$ denotes the outward unit normal of $T\in\mathcal T_h^3$. We write $\mathcal P^{\mathrm t}_{r-2}(f)$ for the space of tangent polynomial vector fields on $f$ of degree at most $r-2$. Lastly, for vector proxies of $2$-forms, letting $b_T$ denote the barycenter of $T\in\mathcal T_h^3$, and denoting $\tilde{\lambda}_k = \lambda_k - \lambda_k(b_T)$, the operator $K_{T}$ is determined by
\begin{equation}
\label{eq:cone3d--koszul}
K_{T}(\tilde{\lambda}^\alpha \nabla\lambda_i\times \nabla\lambda_j)=\frac{\tilde{\lambda}^\alpha}{\abs{\alpha}+2} (\nabla\lambda_i\times\nabla\lambda_j)\times (x-b_T).
\end{equation}

\begin{corollary}[Curl potentials from cones]
\label{cor:cone3d--P2}
Let $\mathbf u\in\mathsf{RT}_r(\mathcal T_h)$. Then $\mathscr P^2_{r,\mathrm{co}}\mathbf u$ is the unique element of $\mathsf{N1}_r(\mathcal T_h)$ such that, for all $e\in\mathcal T_h^1$, $f\in\mathcal T_h^2$ and $T\in\mathcal T_h^3$,
\begin{subequations}
\label{eq:cone3d--P2-def}
\begin{alignat}{2}
&\int_e(\mathscr P^2_{r,\mathrm{co}}\mathbf u\cdot\mathbf t_e)\,q
=\int_{\mathbf C(e)}(\mathbf u\cdot\boldsymbol \nu_{\mathbf C(e)})\,\ext_eq
&\qquad&\forall q\in\mathcal P_{r-1}(e),
\label{eq:cone3d--P2-edge}\\
&\int_f(\mathscr P^2_{r,\mathrm{co}}\mathbf u\times\mathbf q)\cdot\boldsymbol\nu_f
=\int_{\mathbf C(f)}\mathbf u\cdot \ext_f\mathbf q
&&\forall\mathbf q\in\mathcal P^{\mathrm t}_{r-2}(f),
\label{eq:cone3d--P2-face}\\
&\int_T\mathscr P^2_{r,\mathrm{co}} \mathbf u\cdot\mathbf q
=\int_T\mathbf u\cdot K_T\mathbf q-\sum_{f\triangleleft T}(\boldsymbol \nu_T\cdot\boldsymbol\nu_f)\int_{\mathbf C(f)}\mathbf u\cdot\ext_f(K_T\mathbf q)
&&\forall\mathbf q\in[\mathcal P_{r-3}(T)]^3 .
\label{eq:cone3d--P2-cell}
\end{alignat}
\end{subequations}
\end{corollary}

Note that $\boldsymbol{\nu}_T\cdot \boldsymbol{\nu}_f = \pm 1$. For $r=1$ the only constraint is \eqref{eq:cone3d--P2-edge} with $q$ constant, and we recover the lowest order operator \eqref{eq:conext--poincare-lowest-order}.

\subsection{Setting of the general construction}\label{section:conext--general-setting}

We now discuss how the construction generalizes past the setting of star-shaped domains. By re-parametrizing the cone $\mathbf C(\sigma)$ we will obtain simplified formulas for the right-hand side of \eqref{eq:conext--poincare-def-subeq}. This new perspective also makes clear how the cone Poincar\'e operator \eqref{eq:conext--poincare-def-subeq} relates to the classic Poincar\'e operator \eqref{eq:prelim--classical-poincare}. Throughout the rest of this section, we assume $\Omega$ is Lipschitz contractible in the following sense.
\begin{definition}[Lipschitz contractible domain]
\label{def:conext--Lipschitz-contraction}
    We say a domain $\Omega\subset\mathbb R^n$ is \emph{Lipschitz contractible} to a point $p_0\in\Omega$ if there exists a Lipschitz map $\Phi: [0,1]\times \overline \Omega \to \overline \Omega$ such that $\Phi(1,x) = x$ and $\Phi(0,x) = p_0$, and such that $\Phi(t, \cdot\,)$ is locally bi-Lipschitz for all $t \in (0, 1]$. 
\end{definition}
Let now $\Phi$ be a fixed Lipschitz contraction of $\Omega$ to a point $p_0$. Using $\Phi$, we will define a discrete Poincar\'e operator which is a generalization of \cref{def:conext--poincare}.

\begin{definition}[Extrusion of a simplex]
\label{def:conext--cone}
For $\sigma\in\mathcal T_h^\ell$, we define the \emph{extrusion} of $\sigma$ to be the ``prism''
\begin{equation}
\label{eq:conext--cone-of-simplex}
\cone(\sigma)\coloneqq[0,1]\times\sigma ,
\end{equation}
oriented so that $dt\wedge d\lambda_1\wedge\cdots\wedge d\lambda_\ell$ is positive. 
\end{definition}

\begin{definition}
Let $i_{\cone(\sigma)}:\cone(\sigma)\hookrightarrow [0,1]\times \overline\Omega$ be the inclusion map. We define
\begin{equation}
    \Phi_\sigma \coloneqq \Phi\circ i_{\cone(\sigma)}: \cone(\sigma)\to \overline \Omega.
\end{equation}
Note that $\Phi_\sigma$ collapses the face $t=0$ to $p_0$. We also set $\mathbf{C}(\sigma):=\im \Phi_\sigma$, which generalizes the previous notion of a geometric cone.

\end{definition}

\begin{example}[Model case: star-shaped domains]
\label{ex:conext--star-shaped}
Let $\Omega$ be star-shaped with respect to $p_0$. Then the \emph{straight-line contraction}
\begin{equation}
\label{eq:conext--cone-map}
\Phi(t,x)\coloneqq(1-t)p_0+tx
\end{equation}
is a Lipschitz contraction of $\Omega$ onto $p_0$ in the sense of \cref{def:conext--Lipschitz-contraction}. 
\end{example}

The natural idea is to extend forms which live on $\sigma$ directly to $\mathbf{C}(\sigma)$. However, it turns out to be simpler to work with forms on $\cone(\sigma)$ instead. Indeed, this lets us provide a unified treatment for all cell dimensions. When $\Phi_\sigma$ is injective on $(0,1]\times\sigma$, as in \cref{ex:conext--star-shaped} for $\dim\sigma<n$ and $p_0$ affinely independent from $\sigma$, extensions to $\cone(\sigma)$ naturally correspond to extensions to $\mathbf{C}(\sigma)$.

This idea also relates to the classical Poincar\'e operator. In computing the canonical interpolant of $\mathscr P_{p_0}^{k+1}u$, we have integrals of the form
\begin{equation}
\label{eq:conext--interpolant-idea}
    \int_\sigma \mathscr P_{p_0}^{k+1} u\wedge q=\int_\sigma\bigl(\int_0^1\cdots\quad\bigr) \wedge q.
\end{equation}
By ``pulling'' $q$ into the inner-most integral, we can write this as an integral over $[0,1]\times \sigma = \cone(\sigma)$. In particular, our discrete Poincar\'e operator in the general setting is identical to the canonical interpolant of $\mathscr P_{p_0}^{k+1}u$ for smooth $u$ (see \cref{section:conext--classical-relation}).

\subsection{General discrete Poincar\'e operator}\label{section:conext--general-construction}

For each $\sigma \in\mathcal T_h$ let $\pr_\sigma: \cone(\sigma)\to \sigma$ denote the projection onto the spatial component $(t,y)\mapsto y$. The \emph{prism extension} on $\sigma$ is the pullback $\pr^*_\sigma:H\Lambda^k(\sigma)\to H\Lambda^k(\cone(\sigma))$. This may be thought of as a simple extension by constant. With the prism extension in hand, we can now define the discrete Poincar\'e operator in the general setting.

\begin{definition}[General discrete Poincar\'e operator]
\label{def:conext--poincare-general}
For $0\leq k\leq n-1$, the map $\mathscr P_{r, \Phi}^{k+1}:\mathcal P_r^-\Lambda^{k+1}(\mathcal T_h)\to \mathcal P_r^-\Lambda^k(\mathcal T_h)$ is defined by prescribing its canonical dof as
\begin{equation}
\label{eq:conext--poincare-def-general}
\int_\sigma \mathscr P_{r, \Phi}^{k+1} u\wedge q = \int_{\cone(\sigma)} \Phi_\sigma^*u\wedge \pr^*_\sigma q \qquad \forall q\in \mathcal P_{r+k-\ell-1}\Lambda^{\ell-k}(\sigma),
\end{equation}
for all $\sigma\in \mathcal T_h^\ell$ and all $\ell\geq k$.
\end{definition}

If $\Phi_\sigma$ is injective on $(0,1]\times\sigma$, then a change of variables yields
\begin{equation}\label{eq:conext--prism-cone-relation}
    \int_{\cone(\sigma)} \Phi_\sigma^*u\wedge \pr^*_\sigma q = \int_{\mathbf C(\sigma)} u \wedge (\Phi_\sigma^{-1})^* \pr_\sigma^* q.
\end{equation}
When $\Phi$ is the straight-line contraction, this reproduces \eqref{eq:conext--poincare-def} with the rational extension, and so \cref{def:conext--poincare-general} is indeed a generalization of \cref{def:conext--poincare} as long as the top cell dofs agree (see \cref{rem:conext--cone-Phi-same}). One can view the equality \eqref{eq:conext--prism-cone-relation} as a Duffy transform \cite{Duffy:82} mapping the trace of $u$ on $\mathbf{C}(\sigma)$ to the prism $\mathbf{P}(\sigma)$.

We now want to show that the operator of \cref{def:conext--poincare-general} satisfies the homotopy identity. As in \cref{section:conext--cone-star-shaped}, the key is an integration by parts result. For that we need to assume some regularity. We say that $u\in H\Lambda^k(\Omega)$ is bounded near $p_0$ if there exists an open neighbourhood $U$ of $p_0$ such that $\tr_U u \in \mathrsfs W^k_{\infty,\infty}(U)$.

\begin{lemma}
\label{lem:conext--regularity-assumption}
If $u\in H\Lambda^k(\Omega)$ is bounded near $p_0$, then $\Phi^*u \in H\Lambda^k([0, 1]\times \Omega)$.
    
\end{lemma}

\begin{proof}

Since $d$ commutes with pullbacks, it suffices to show that $\Phi^*u$ and $\Phi^*du$ are $L^2$-integrable. Let $U$ be an open neighbourhood of $p_0$ such that $\tr_U u \in \mathrsfs W^k_{\infty,\infty}(U)$. Since $\Phi$ is globally Lipschitz there exists $\delta > 0$ such that $\Phi(t, y) \in U$ holds for all $(t,y) \in [0, \delta]\times \Omega$. Writing out $u$ in terms of the standard basis for $\Lambda^k(\RR^{n,*})$ yields
\begin{equation}
    \normdelim{\Phi^*u}_{L^\infty([0,\delta]\times \Omega)} \leq \normdelim{D\Phi}_{L^\infty([0,\delta]\times \Omega)}^k \normdelim{u}_{L^\infty(U)},
\end{equation}
and likewise for $\normdelim{\Phi^*du}_{L^\infty([0,\delta]\times\Omega)}$. Since $\normdelim{D\Phi}_\infty$ is bounded by the Lipschitz constant for $\Phi$ this gives us $\Phi^*u\in \mathrsfs W^k_{\infty, \infty}([0,\delta]\times \Omega)$, and because $\Omega$ is bounded $\mathrsfs W^k_{\infty, \infty}([0,\delta]\times \Omega)\subset H\Lambda^k([0,\delta]\times \Omega)$. As for the remaining bit, $\Phi^*u$ and $\Phi^*du$ are $L^2$-integrable on $[\delta,1]\times \Omega$ due to \cite[Theorem 2.2]{Goldshtein1982}, using that $\Phi^*$ is locally bi-Lipschitz on the quadrant $[\delta, 1]$.
\end{proof}

\begin{lemma}[Integration by parts over the extrusion]
\label{lem:conext--parts}
Let $\sigma\in\mathcal T_h^\ell$, $q\in H\Lambda^{\ell-k}(\sigma)$, and let $u\in H\Lambda^k(\Omega)$ bounded near $p_0$. For $k\geq 1$ it holds
\begin{subequations}
\label{eq:conext--parts}
\begin{equation}
\label{eq:conext--parts-k}
\int_{\cone(\sigma)}\!\bigl(\Phi_\sigma^*du\wedge \pr^*_\sigma q
+(-1)^k\Phi_\sigma^*u\wedge \pr^*_\sigma dq\bigr) +\sum_{\tau\triangleleft\sigma}o_{\sigma\tau}\int_{\cone(\tau)}\Phi_\tau^*u\wedge\pr^*_\tau\tr_\tau q =\int_\sigma u\wedge q .
\end{equation}
For $k=0$ we instead have
\begin{equation}
\label{eq:conext--parts-0}
\int_{\cone(\sigma)}\!\Phi_\sigma^*du\wedge \pr^*_\sigma q
=\int_\sigma u\wedge q-u(p_0)\!\int_\sigma q .
\end{equation}
\end{subequations}
\end{lemma}

\begin{proof}

Write $\sigma_1\coloneqq\{1\}\times\sigma$ and $\sigma_0\coloneqq\{0\}\times\sigma$, both identified with $\sigma$ through $(t,y)\mapsto y$. With the orientations induced on its faces, the boundary of $\cone(\sigma)$ is
\begin{equation}
\label{eq:conext--prism-boundary}
\partial\cone(\sigma)=\sigma_1-\sigma_0-\sum_{\tau\triangleleft\sigma}o_{\sigma\tau}\cone(\tau),
\end{equation}
which can be seen, for instance, by partitioning $\cone(\sigma)$ into simplices as in the proof of \cite[Theorem 2.10]{Hatcher:02}. Let $\eta\coloneqq\Phi_\sigma^*u\wedge\pr^*_\sigma q$. 
Since pointwise $\ipdelim{\pr_\sigma^*q, \pr_\sigma^*q} = \ipdelim{q,q}$, it follows from Fubini's theorem that $\pr_\sigma^*$ maps $L^p\Lambda^k(\sigma) \to L^p\Lambda^k(\cone(\sigma))$, and thus $\pr^*_\sigma q\in H\Lambda^{\ell-k}(\cone(\sigma))$. We have that $\Phi^*_\sigma u\in H\Lambda^k(\cone(\sigma))$ by \cref{lem:conext--regularity-assumption}, so put together this gives us $\eta\in \mathrsfs W_{1,1}^{\ell}(\cone(\sigma))$. From here we simply want to apply Stokes' theorem to $\eta$ in order to obtain \eqref{eq:conext--parts}. In light of \eqref{eq:conext--prism-boundary}, it states that
\begin{equation}
\label{eq:conext--parts-stokes}
    \int_{\cone(\sigma)} d\eta = \int_{\sigma_1} \tr_{\sigma_1}\eta - \int_{\sigma_0} \tr_{\sigma_0} \eta - \sum_{\tau\triangleleft \sigma} o_{\sigma\tau} \int_{\cone(\tau)} \tr_{\cone(\tau)} \eta.
\end{equation}
We proceed by individually expanding each term.

Using that $d$ commutes with pullbacks, the Leibniz rule yields 
\begin{equation}
\label{eq:conext--leibniz}
d\eta=\Phi_\sigma^*du\wedge\pr^*_\sigma q+(-1)^k\Phi_\sigma^*u\wedge\pr^*_\sigma dq,
\end{equation}
which is the leftmost integrand in \eqref{eq:conext--parts} (if $k=0$ then $dq=0$). For the boundary term we find that $\tr_{\cone(\tau)}\Phi^*_\sigma = \Phi^*_\tau$ by unravelling the definition $i_{\cone(\tau)}^*\circ(i_{\cone(\sigma)}^*\circ\Phi^*)$, and similarly $\tr_{\cone(\tau)}\pr^*_\sigma = \pr_\tau^*\circ \tr_\tau$ by the fact $\pr_\sigma\circ i_{\cone(\tau)\hookrightarrow\cone(\sigma)} = i_{\tau\hookrightarrow \sigma}\circ \pr_\tau$. Taken together, it holds
\begin{equation}
\tr_{\cone(\tau)}\eta = \Phi^*_\tau u\wedge \pr_\tau^* \tr_\tau q,
\end{equation}
and this is the second left-hand side integrand of \eqref{eq:conext--parts-k} (if $k=0$ then $\tr_\tau q =0$). It remains to compute the contributions on $\sigma_1$ and $\sigma_0$. For the former, we use that both $\Phi_\sigma$ and $\pr_\sigma$ act on $\sigma_1$ as the identity, so that $\tr_{\sigma_1}\eta \equiv u\wedge q$. For the latter, since $\Phi_\sigma$ is constant on $\sigma_0$ we only see a contribution when $k=0$:
\begin{equation}
    \tr_{\sigma_0}\eta \equiv \begin{cases}
        u(p_0)\wedge q & \text{if $k=0$}, \\
        0 & \text{if $k\geq 1$}.
    \end{cases}
\end{equation}
Plugging everything back into \eqref{eq:conext--parts-stokes} and rearranging yields \eqref{eq:conext--parts}. The point evaluation $u(p_0)$ is well-defined for $k=0$ because $u$ is Lipschitz at $p_0$ in this case.
\end{proof}

\begin{corollary}[Homotopy formula in the general setting]
    The statement of \cref{thm:conext--homotopy-formula} holds for the operator of \cref{def:conext--poincare-general}.
\end{corollary}
\begin{proof}
    By the same line of arguments as in the proof of \cref{thm:conext--homotopy-formula}, using \cref{lem:conext--parts} in place of \cref{lem:conext--cone-stokes}. In particular, we do not need a separate treatment for the $k=n$ case.
\end{proof}

\subsection{Relationship between the discrete and continuous Poincar\'e operator}\label{section:conext--classical-relation}

We now show how the discrete Poincar\'e operator of \cref{def:conext--poincare-general} can be identified with the canonical interpolant of the classical Poincar\'e operator. The standard Poincar\'e operator for a smoothly contractible domain is given by
\begin{equation}
    \mathscr P_{p_0}^{k+1}u = \int_0^1 \boldsymbol{\partial}_t \lrcorner\, (\Phi^* u)\,dt,
\end{equation}
which is to be interpreted pointwise, with the integrand taking  values in $\Lambda^k(T_x^* \Omega)$ (see \cite[Lemma 17.9]{Lee:13}). We show that $\mathscr P_{r,\Phi}^{k+1} = I^k_r\circ\mathscr P_{p_0}^{k+1}$, with $I^k_r:\Lambda^k(\Omega)\to \mathcal P_r^-\Lambda^k(\mathcal T_h)$ denoting the canonical interpolant. Since $d$ commutes with $I^k_r$, and since $\mathscr P_{p_0}^\bullet$ satisfies the homotopy formula, this relation implies the homotopy identity for $\mathscr P_{r, \Phi}^{k+1}$ in the smooth setting. 

\begin{proposition}\label{prop:conext--classical-interpolant}
    For all $\sigma \in \mathcal T_h$, $u\in \Lambda^{k+1}(\Omega)$, and $q\in \Lambda^{\dim\sigma-k}(\sigma)$, it holds
\begin{equation}
\label{eq:conext--interp-identity}
\int_\sigma\mathscr P_{p_0}^{k+1}u\wedge q 
= \int_{\cone(\sigma)}\Phi_\sigma^*u\wedge\pr^*_\sigma q.
\end{equation}
\end{proposition}

\begin{proof}
    To show \eqref{eq:conext--interp-identity}, we compute for $\sigma\in \mathcal T_h$
\begin{align}
\int_\sigma\mathscr P_{p_0}^{k+1}u\wedge q 
&= \int_\sigma\bigl(\int_0^1\boldsymbol{\partial}_t\lrcorner\, \Phi_\sigma^*u\,dt\bigr) \wedge q \\ \nonumber
&= \int_\sigma\bigl(\int_0^1(\boldsymbol{\partial}_t\lrcorner\, \Phi_\sigma^*u)\wedge \pr^*_\sigma q\,dt\bigr) \\ \nonumber
&= \int_\sigma\bigl(\int_0^1\boldsymbol{\partial}_t\lrcorner\,( \Phi_\sigma^*u\wedge \pr^*_\sigma q) \,dt \bigr)
+ (-1)^{k}\int_\sigma\bigl(\int_0^1 \Phi_\sigma^*u\wedge (\underbrace{\boldsymbol{\partial}_t\lrcorner\, \pr^*_\sigma q}_{=\,0}) \,dt\bigr) \\ \nonumber
&= \int_{[0,1]\times \sigma} \Phi_\sigma^*u\wedge \pr^*_\sigma q.
\end{align}
The final equality is Fubini's theorem, using that $\boldsymbol{\partial}_t \lrcorner\, (dt\wedge dy^1\wedge \cdots dy^\ell) = dy^1\wedge \cdots \wedge dy^\ell$.
\end{proof}

\begin{remark}\label{rem:conext--cone-Phi-same}
In the model case of \cref{ex:conext--star-shaped}, and under the assumption that $p_0$ is in general position, it holds $\mathscr P_{r, \mathrm{co}} = \mathscr P_{r, \Phi} = I_r\circ\mathscr P_{p_0}$ provided $K_\sigma$ is the scaled Koszul operator $K_{p_0}$. Indeed, using the fact that $dK_{p_0}+K_{p_0}d=\id$, we get that
\begin{equation}
    d(\Phi^*_\sigma u\wedge \pr^*_\sigma K_{p_0} q) 
    = \Phi^*_\sigma du\wedge \pr^*_\sigma K_{p_0}q + (-1)^{k+1}\Phi^*_\sigma u\wedge \pr^*_\sigma q + (-1)^k\Phi^*_\sigma u\wedge \pr^*_\sigma K_{p_0}dq.
\end{equation}
One checks that $\kappa_{p_0} q=0\implies \Phi^*_\sigma u\wedge \pr_\sigma^* q=0$, and consequently both the first and the latter term on the right-hand side vanish, so that only the middle term is non-zero. Then the equality \eqref{eq:conext--poincare-def-general} for $\sigma\in\mathcal T_h^n$ follows from Stokes' theorem by the same line of arguments as in the proof of \cref{lem:conext--parts}.
\end{remark}

\section{Discrete Bogovski\u{\i} operators}
\label{section:bogovskii}

We now turn to the complex with homogeneous boundary conditions,
\begin{equation}
\label{eq:bogovskii--ring-complex}
  \begin{tikzcd}[column sep = large, ampersand replacement=\&]
    \mathring{\mathcal P}_r^-\Lambda^0(\mathcal T_h) \arrow[r, "d"] \&
    \mathring{\mathcal P}_r^-\Lambda^1(\mathcal T_h) \arrow[r, "d"] \&
    \cdots \arrow[r, "d"] \&
    \mathring{\mathcal P}_r^-\Lambda^n(\mathcal T_h),
  \end{tikzcd}
\end{equation}
and seek homotopy operators
$\mathscr B^k:\mathring{\mathcal P}_r^-\Lambda^k(\mathcal T_h)\to \mathring{\mathcal P}_r^-\Lambda^{k-1}(\mathcal T_h)$
satisfying
\begin{equation}
\label{eq:bogovskii--homotopy-formula}
d\mathscr B^ku+\mathscr B^{k+1}du = u - \rho^k u .
\end{equation}
These are called discrete Bogovski\u{\i} operators, in reference to \cite{Bogovskii:79}. The cohomology of \eqref{eq:bogovskii--ring-complex} is the relative de Rham cohomology of $(\overline\Omega,\partial\Omega)$ \cite[Section~6]{Arnold-Falk-Winther:10}, which is isomorphic to $H^{n-k}(\Omega)$ by the Poincar\'e-Lefschetz duality \cite[Theorem 3.43]{Hatcher:02}. We can therefore think of $\rho^\bullet$ as playing the role of $\pi^\bullet$, but with $\im\rho^k\cong H^k(\overline\Omega,\partial\Omega)\cong H^{n-k}(\Omega)$. In particular $\rho^n \neq 0$.

Both constructions of the previous two sections have Bogovski\u{\i} analogues.

\subsection{Discrete Bogovski\u{\i} operator from the decomposition}
\label{section:bogovskii--decomposition}

The operator $B^\bullet$ already preserves homogeneous boundary conditions, so only the lowest order factor of the decomposition \eqref{eq:lc--relative-decomposition} needs to be replaced.

\begin{proposition}
\label{prop:bogovskii--decomposition}
Let $\mathscr B_{\mathrm W}^k:\mathring{\mathcal P}_1^-\Lambda^k(\mathcal T_h)\to\mathring{\mathcal P}_1^-\Lambda^{k-1}(\mathcal T_h)$ be a homotopy operator on the Whitney complex with homogeneous boundary conditions, with cohomology projection $\rho_{\mathrm W}^\bullet$. Then
\begin{equation}
\label{eq:bogovskii--decomposition-def}
\mathscr B_{r,\mathrm{dec}}^k u \coloneqq \mathscr B^k_{\mathrm W} I^ku + B^k(\id-I^k)u
\end{equation}
defines a homotopy operator on \eqref{eq:bogovskii--ring-complex} satisfying \eqref{eq:bogovskii--homotopy-formula} with $\rho^k=\rho^k_{\mathrm W}I^k$.
\end{proposition}

\begin{proof}
By \cref{lem:lc--interpolant}\,(iv), $I^k$ maps $\mathring{\mathcal P}_r^-\Lambda^k(\mathcal T_h)$ into $\mathring{\mathcal P}_1^-\Lambda^k(\mathcal T_h)$, so $\id-I^k$ maps it into $\mathring{\mathcal M}_r^-\Lambda^k(\mathcal T_h)$. That space is preserved by $B^k$, and \eqref{eq:lc--ag-homotopy} holds because $\mathring{\mathcal M}_r^-\Lambda^k(\mathcal T_h)\subset\mathcal M_r^-\Lambda^k(\mathcal T_h)$. The computation in the proof of \cref{prop:lc--abstract-construction} now applies verbatim, with \eqref{eq:lc--relative-decomposition} in place of \eqref{eq:lc--wr-decomposition-global}.
\end{proof}

One choice of $\mathscr B^\bullet_{\mathrm W}$ is the operator $\mathscr H^\bullet_V$ of \cref{prop:prelim--abstract-bstar} for the complex \eqref{eq:bogovskii--ring-complex} at $r=1$ with the $L^2(\Omega)$ inner product. That is, the solution operator of \eqref{eq:prelim--hodge-laplace} as posed over $\mathring{\mathcal P}_1^-\Lambda^\bullet(\mathcal T_h)$. Its cohomology projection $\rho^\bullet_{\mathrm W}$ is then the $L^2$-orthogonal projection onto the harmonic forms $\mathfrak H\mathring{\mathcal P}_1^-\Lambda^\bullet(\mathcal T_h)$ of \eqref{eq:bogovskii--ring-complex}. Alternatives are discussed in \cite{guzman-potu:26}.

\subsection{Discrete Bogovski\u{\i} operator from the cone construction}
\label{section:bogovskii--cone}

Throughout this subsection we place ourselves in the model case of \cref{ex:conext--star-shaped}: we suppose $\Omega$ is star-shaped with respect to $p_0$ and $\Phi$ is the straight-line contraction \eqref{eq:conext--cone-map}. In this setting the cone construction of \cref{section:cone-extension} has a counterpart in which the cone $\mathbf {C}(\sigma)$, swept out by contracting $\sigma$ towards $p_0$, is replaced by the ``shadow'' of $\sigma$ as seen from $p_0$, swept out by pushing $\sigma$ away from $p_0$ until it leaves $\Omega$; the lowest order case is the construction of \cite[Section 7]{guzman-potu:26}. 

In addition to star-shapedness we assume that $p_0$ lies in the interior of a cell $T_*\in\mathcal T_h^n$, whose increasing vertex enumeration we assume, relabeling if necessary, to orient it positively with respect to $\RR^n$. These are the same assumption as are made in \cite[Section 7]{guzman-potu:26}. Because of this star-shapedness, if $x\in\partial\Omega$ and $t>1$, then $p_0+t(x-p_0)\notin\Omega$. 

\begin{definition}[Shadow of a simplex]
\label{def:bogovskii--shadow}
For $\sigma\in\mathcal T_h^\ell$, we define the \emph{shadow} of $\sigma$ to be the half-infinite prism
\begin{equation}
\label{eq:bogovskii--shadow-of-simplex}
\cone^\vee(\sigma)\coloneqq[1,\infty)\times\sigma ,
\end{equation}
oriented so that $dt\wedge d\lambda_1\wedge\cdots\wedge d\lambda_\ell$ is positive. 
\end{definition}

\begin{definition}
Let $i_{\cone^\vee(\sigma)}:\cone^\vee(\sigma)\hookrightarrow [1, \infty)\times\RR^n$ be the inclusion. The contraction $\Phi(t,y)= (1-t) p_0 + ty$ extends to a map $\Phi^\vee:[1, \infty)\times \RR^n \to \RR^n$ in the obvious way, and we write 
\begin{equation}
    \Phi^\vee_\sigma \coloneqq \Phi^\vee \circ i_{\cone^\vee(\sigma)} : \cone^\vee(\sigma)\to \RR^n.
\end{equation}
The map $\Phi^\vee_\sigma$ is the continuation of $\Phi_\sigma$ beyond $t=1$, and its image $\Phi^\vee_\sigma(\cone^\vee(\sigma)) \cap\overline\Omega$ is the geometric shadow of $\sigma$.
\end{definition}

As before, we write $\pr_\sigma$ to mean the map $(t,y)\mapsto y$, except that here the domain is $\cone^\vee(\sigma)$ instead of $\cone(\sigma)$. For a form $u$ on $\Omega$, we write $\tilde u$ to mean its extension by zero to $\RR^n$. We present the integration by parts formula over the shadow, the counterpart of \cref{lem:conext--parts}, from which the construction follows exactly as in \cref{section:conext--general-construction}.

\begin{proposition}[Integration by parts over the shadow]
\label{lem:bogovskii--shadow-parts}
Let $\sigma\in\mathcal T_h^\ell$, $q\in \Lambda^{\ell-k}(\sigma)$, and $u\in\mathring{\mathcal P}_r^-\Lambda^k(\mathcal T_h)$. 
If $\sigma\neq T_*$ or $k\leq n-1$, then
\begin{subequations}
\label{eq:bogovskii--shadow-parts}
\begin{equation}
\label{eq:bogovskii--shadow-parts-k}
\int_{\cone^\vee(\sigma)}\!\!\bigl((\Phi^\vee_\sigma)^*\widetilde{du}\wedge\pr^*_\sigma q
+(-1)^k(\Phi^\vee_\sigma)^*\tilde u\wedge\pr^*_\sigma dq\bigr)
+\sum_{\tau\triangleleft \sigma} o_{\sigma\tau}\int_{\cone^\vee(\tau)}\!\!(\Phi^\vee_{\tau})^*\tilde u\wedge\pr^*_{\tau}\tr_\tau q
=-\int_\sigma u\wedge q.
\end{equation}
If $\sigma = T_*$ and $k=n$ we pick up an additional term on the right-hand side
\begin{equation}
\label{eq:bogovskii--shadow-parts-top}
(-1)^n\!\!\int_{\cone^\vee(T_*)}\!\!(\Phi^\vee_{T_*})^*\tilde u\wedge\pr^*_{T_*} dq
+\sum_{\tau\triangleleft T_*} o_{T_{*}\tau}\int_{\cone^\vee(\tau)}\!\!(\Phi^\vee_{\tau})^*\tilde u\wedge\pr^*_{\tau}\tr_\tau q
= -\int_{T_*} u\wedge q + q(p_0)\!\int_\Omega u .
\end{equation}
\end{subequations}
\end{proposition}

\begin{proof}

Since the shadow is unbounded, Stokes' theorem will be applied the truncated shadow $[1, R]\times\sigma$ for sufficiently large $R$. We write $\sigma_1\coloneqq\{1\}\times\sigma$ and $\sigma_R \coloneqq \{R\}\times \sigma$. By the same argument as before, it holds
\begin{equation}
\label{eq:bogovskii--shadow-boundary}
\partial\bigl([1, R]\times\sigma\bigr)
= \sigma_R-\sigma_1-\sum_{\tau \triangleleft \sigma}o_{\sigma\tau} \bigl([1,R]\times \tau\bigr).
\end{equation}
With $\eta\coloneqq(\Phi^\vee_\sigma)^*\tilde u\wedge\operatorname{pr}_\sigma^*q$, Stokes' theorem yields
\begin{equation}
    \int_{[1, R] \times \sigma} d \eta = \int_{\sigma_R} \tr_{\sigma_R}\eta - \int_{\sigma_1}\tr_{\sigma_1} \eta - \sum_{\tau\triangleleft \sigma} o_{\sigma\tau} \int_{[1,R]\times \tau}\tr_{\cone^\vee(\tau)}\eta.
\end{equation}
We proceed by computing each term.

Leibniz rule and commutativity of $d$ with the pullback yields
\begin{equation}
\label{eq:bogovskii--shadow-parts-Leibniz}
    d\eta = (\Phi_\sigma^\vee)^*d\tilde u \wedge \pr^*_\sigma q + (-1)^k (\Phi^\vee_\sigma)^*\tilde u\wedge \pr^*_\sigma dq.
\end{equation}
Since the trace of $u$ vanishes on $\partial \Omega$ we have that $d\tilde u = \widetilde{du}$, and thus \eqref{eq:bogovskii--shadow-parts-Leibniz} corresponds to the left-hand side of \eqref{eq:bogovskii--shadow-parts} for sufficiently large $R$ (If $k=n$ then $du = 0$). As for the boundary terms, 
\begin{equation}
    \tr_{\cone^\vee(\tau)}\eta = (\Phi^\vee_\tau)^*\tilde u\wedge \pr^*_\tau \tr_\tau q
\end{equation}
holds by the same argument as before, as well as $\tr_{\sigma_1} \eta\equiv u\wedge q$. All that remains is to compute $\tr_{\sigma_R}\eta$.

Using the formula for the pullback of a composition we have that $\tr_{\sigma_R}(\Phi^\vee_\sigma)^*u = D_R^*u$, where $D_R\coloneqq \Phi^\vee_\sigma \circ i_{\sigma_R}$. It admits the explicit formula $D_R(y)=(1-R)p_0 + Ry$, with inverse $D^{-1}_R(x)= p_0 + R^{-1}(x-p_0)$. Thus we compute $(D^{-1}_R)^*dy=R^{-1}dx$, and in particular $(D^{-1}_R)^*\tr_{\sigma_R}\pr^*_\sigma q = R^{k-\ell}q\circ D_R^{-1}$ (recall that $\tr_{\sigma_R}\pr^*_\sigma$ is effectively the identity). Using the change of variables $y=D_R^{-1}(x)$ yields
\label{eq:bogovskii--shadow-parts-limit-integral}
\begin{align}
    \int_{\sigma_R}\tr_{\sigma_R} \eta 
    &= \int_{\sigma_R} \tr_{\sigma_R}(\Phi^\vee_R)^*\tilde u \wedge \tr_{\sigma_R}\pr^*_\sigma q \\ \nonumber
    &= \int_{D_R(\sigma)} (D_R^{-1})^*(D^*_R\tilde u \wedge \tr_{\sigma_R}\pr^*_\sigma q). \\
    \nonumber
    &= \int_{D_R(\sigma)} \tilde u \wedge R^{k-\ell} q\circ D_R^{-1}.
\end{align}
If $\sigma\ne T_*$ then for sufficiently large $R$ the image of $D_R$ lies entirely outside of $\Omega$, and the integral vanishes. Consequently this term does not appear in \eqref{eq:bogovskii--shadow-parts-k}. Suppose now that $\sigma=T_*$, in which case $\ell=n$ and $q\in \Lambda^{n-k}(T_*)$. As $\tilde u$ has compact support and $q$ is bounded, for $k\leq n-1$ the integral tends to zero in the limit $R\to \infty$. Otherwise, if $k=n$ then the integral converges to $\int_{\RR^n}\tilde u\wedge q(p_0)$ by dominated convergence, which is the additional term in \eqref{eq:bogovskii--shadow-parts-top}. 
\end{proof}

\begin{definition}[Bogovski\u{i} operator]
\label{def:bogovskii--cone-operator}
For $0\leq k\leq n-1$ define $\mathscr B_{r,\mathrm{co}}^{k+1}:\mathring{\mathcal P}_r^-\Lambda^{k+1}(\mathcal T_h)\to\mathring{\mathcal P}_r^-\Lambda^k(\mathcal T_h)$ by
\begin{equation}
\label{eq:bogovskii--cone-def}
\int_\sigma \mathscr B_{r,\mathrm{co}}^{k+1} u\wedge q \coloneqq -\int_{\cone^\vee(\sigma)} (\Phi^\vee_\sigma)^*\tilde u\wedge \pr^*_\sigma q \qquad \forall q\in \mathcal P_{r+k-\ell-1}\Lambda^{\ell-k}(\sigma),
\end{equation}
for all $\sigma\in \mathcal T_h^\ell$ and all $\ell\geq k$.
\end{definition}

One sees that \eqref{eq:bogovskii--cone-def} maps to the space with homogeneous boundary conditions. Indeed, if $\sigma\subset\partial\Omega$, then for every $x\in\sigma$ and $t>1$ the point $\Phi^\vee_\sigma(t,x)=p_0+t(x-p_0)$ does not lie in $\Omega$, by the star-shapedness consequence noted at the start of this subsection. Hence $(\Phi^\vee_\sigma)^*\tilde u$ vanishes almost everywhere on $[1,\infty)\times\sigma$, and the integral on the right of \eqref{eq:bogovskii--cone-def} is zero. All dof attached to faces in $\partial\mathcal T_h$ are thus zero, which forces $\tr_{\partial\Omega}\mathscr B^{k+1}_{r,\mathrm{co}}u=0$.

\begin{theorem}[Bogovski\u{\i} operator]
\label{thm:bogovskii--cone}
Let $\rho^k\coloneqq0$ for $k\leq n-1$. For $u\in \mathcal P_r^-\Lambda^n(\mathcal T_h)$ we set $\rho^nu\in\mathring{\mathcal P}^-_r\Lambda^n(\mathcal T_h)$ to be the form vanishing on every cell other than $T_*$, and with the dof
\begin{equation}
\label{eq:bogovskii--rho}
\int_{T_*}\rho^nu\wedge q=q(p_0)\int_\Omega u\qquad\forall q\in\mathcal P_{r-1}\Lambda^0(T_*).
\end{equation}
Then $\rho^\bullet$ is a cohomology projection, and $\mathscr B^\bullet_{r,\mathrm{co}}$ is a homotopy operator associated with it
\begin{equation}
\label{eq:bogovskii--homotopy-formula-cone}
    d\mathscr B^k_{r,\mathrm{co}} u + \mathscr B^{k+1}_{r, \mathrm{co}} du = u - \rho^k u.
\end{equation}
\end{theorem}

\begin{proof}

It is straightforward check that $\rho^\bullet$ is a cohomology projection. The identity \eqref{eq:bogovskii--homotopy-formula} follows by the same argument as in \cref{thm:conext--homotopy-formula}. Fix $\sigma\in\mathcal T_h^\ell$ and an admissible $q$. Stokes' theorem on $\sigma$ together with \eqref{eq:bogovskii--cone-def}, applied exactly as in \eqref{eq:conext--homotopy-proof-term10}--\eqref{eq:conext--homotopy-proof-term21}, expresses $\int_\sigma(d\mathscr B^k_{r,\mathrm{co}}u+\mathscr B^{k+1}_{r,\mathrm{co}}du)\wedge q$ as minus the left hand side of \eqref{eq:bogovskii--shadow-parts-k}, or of \eqref{eq:bogovskii--shadow-parts-top} when $\sigma=T_*$ and $k=n$. For $k=0$, the $dq$- and $\partial\sigma$-terms vanish because $q$ is then of top degree on $\sigma$. \Cref{lem:bogovskii--shadow-parts} therefore gives $\int_\sigma u\wedge q$, except when $\sigma=T_*$ and $k=n$, where it gives $\int_{T_*}u\wedge q-q(p_0)\int_\Omega u$. By unisolvence this implies \eqref{eq:bogovskii--homotopy-formula-cone}.
\end{proof}

\begin{remark}[Geometric realisation of $\mathscr B^{\bullet}_{r,\mathrm{co}}$]
\label{rem:bogovskii--implementation}
The operator $\mathscr B^{k+1}_{r,\mathrm{co}}$ can be interpreted in the same way as $\mathscr P^{k+1}_{r,\mathrm{co}}$ in \cref{section:conext--cone-star-shaped}. Let $p_0$ be in general position with respect to $\mathcal T_h$, and let $\sigma\in\mathcal T_h^\ell$ with $\ell\leq n-1$. Then $\Phi^\vee_\sigma$ is injective on all of $\cone^\vee(\sigma)=[1,\infty)\times\sigma$. The geometric shadow $\mathbf C^\vee(\sigma)\coloneqq\Phi^\vee_\sigma(\cone^\vee(\sigma))\cap\overline\Omega$ is non-degenerate and we can define the extension
\begin{equation}
\label{eq:bogovskii--shadow-ext}
\ext^\vee_\sigma q\coloneqq\bigl((\Phi^\vee_\sigma)^{-1}\bigr)^*\pr^*_\sigma q .
\end{equation}
Unlike in \cref{section:cone-extension}, $\ext^\vee_\sigma q$ is smooth and has no singularity.

Let $\sigma \in\mathcal T_h^\ell$, $u\in \mathring{\mathcal P}_r^-\Lambda^{k+1}(\mathcal T_h)$, and $q\in \mathcal P_{r+k+\ell-1}\Lambda^{\ell-k}(\sigma)$. For $\ell \leq n-1$ it holds
\begin{subequations}
\label{eq:bogovskii--geometric}
\begin{equation}
\label{eq:bogovskii--geometric-low}
    \int_\sigma \mathscr B^{k+1}_{r,\mathrm{co}}u\wedge q =-\int_{\mathbf C^\vee(\sigma)}u\wedge\ext^\vee_\sigma q.
\end{equation}
In the case $\sigma\in\mathcal T_h^n$ we instead get
\begin{equation}
\label{eq:bogovskii--geometric-top}
    \int_\sigma \mathscr B^{k+1}_{r,\mathrm{co}}u\wedge q = (-1)^k\int_\sigma u\wedge K_{p_0}q + (-1)^{k+1}\sum_{\tau\triangleleft \sigma}o_{\sigma\tau}\int_{\mathbf C^\vee(\tau)}u\wedge\ext^\vee_\tau(\tr_\tau K_{p_0}q),
\end{equation}
\end{subequations}
and \eqref{eq:bogovskii--geometric-low}--\eqref{eq:bogovskii--geometric-top} completely determine $\mathscr B^{k+1}_{r\mathrm{co}} u$.

\end{remark}


\section{Implementation in Firedrake}
\label{section:implementation}

We have implemented the local construction of \cref{section:local-construction} in Firedrake \cite{Ham-et-al:23}, for the discretizations of the de Rham complex up to dimension three given by the trimmed and full polynomial forms. Specifically, we have implemented the homotopy operator given by \cref{prop:lc--abstract-construction,prop:lc--global-contraction}, using the strategy of \cref{rem:lc--cost}. To handle the lowest order component we solve the global $r=1$ instance of \eqref{eq:prelim--hodge-laplace}. The implementation supports the following complexes for all admissible choices of polynomial degree:
\[
\begin{array}{c@{\qquad}c}
\text{1D} & \hspace{-21.42em}
\substack{\text{Lagrange}}
\longrightarrow
\substack{\text{Discontinuous Lagrange}}
\\[1.2ex]

\text{2D} & \hspace{-11.23em}
\substack{\text{Lagrange}}
\longrightarrow \hspace{1em}
\substack{\text{Nédélec 1st kind}\\\text{Nédélec 2nd kind}} \hspace{1em}
\longrightarrow
\substack{\text{Discontinuous Lagrange}}
\\[1.8ex]

& \hspace{-11.19em}
\substack{\text{Lagrange}}
\longrightarrow
\substack{\text{Raviart--Thomas}\\\text{Brezzi--Douglas--Marini}}
\longrightarrow
\substack{\text{Discontinuous Lagrange}}
\\[1.8ex]

\text{3D} & \hspace{-1em}
\substack{\text{Lagrange}}
\longrightarrow \hspace{1em}
\substack{\text{Nédélec 1st kind}\\\text{Nédélec 2nd kind}} \hspace{1em}
\longrightarrow
\substack{\text{Raviart--Thomas}\\\text{Brezzi--Douglas--Marini}}
\longrightarrow
\substack{\text{Discontinuous Lagrange}}
\end{array}
\]
All of the code is available at \href{https://github.com/olavtyssvang/discrete-poincare}{github.com/olavtyssvang/discrete-poincare}. For our numerical experiments we only report results for the trimmed spaces.

\subsection{Numerical results}

We have solved the equation $du=f$ over the unit cube $\square^3\subset \mathbb R^3$, for $f$ such that $df=0$. Our point of comparison is the global problem \eqref{eq:prelim--hodge-laplace}, which is implemented using the same mixed formulation as in \eqref{eq:lc--local-problem}. The local algorithm consists of solving the lowest order system \eqref{eq:prelim--hodge-laplace}, as well as applying a matrix of width $\dim \mathcal P_r^-\Lambda^k(\sigma)= O(r^3)$ to the dof on each cell. Therefore, for a mesh consisting of $N$ cells, the higher order correction is of time complexity $O(r^6N)$. 

The $k=1$ instance of \eqref{eq:prelim--hodge-laplace}, as posed over $\mathsf{CG}_r(\mathcal T_h)$ and $\mathsf{N1}_r(\mathcal T_h)$, is solved using the conjugate gradient method with an algebraic multigrid (AMG) preconditioner \cite{falgout-etal:04}. Here we need to take into account the non-zero cohomology $\mathfrak H\mathcal P_r^-\Lambda^0(\square^3) \equiv \RR$. For $k=2$, \eqref{eq:prelim--hodge-laplace} is posed over $\mathsf{N1}_r(\mathcal T_h)$ and $\mathsf{RT}_r(\mathcal T_h)$, and the system is solved using the minimal residual method (MINRES) \cite{choi-etal:11} with a Riesz map preconditioner \cite{kirby-robert:10}. The inverse of the $(\lambda, q)_\Omega$ block is approximated by the Jacobi preconditioner (the inverse of its diagonal), and the inverse of the Riesz block $(dv,dv')_\Omega + (v, v')_\Omega$ is represented by a sparse direct solver (using the LU decomposition). An analogous MINRES/Riesz map strategy is also used for $k=3$, which is posed over $\mathsf{RT}_r(\mathcal T_h)$ and $\mathsf{DG}_{r-1}(\mathcal T_h)$.

\begin{figure}[h]
    \centering
    \begin{subfigure}[t]{0.5\linewidth}
        \centering
        \includegraphics[width=\linewidth]
            {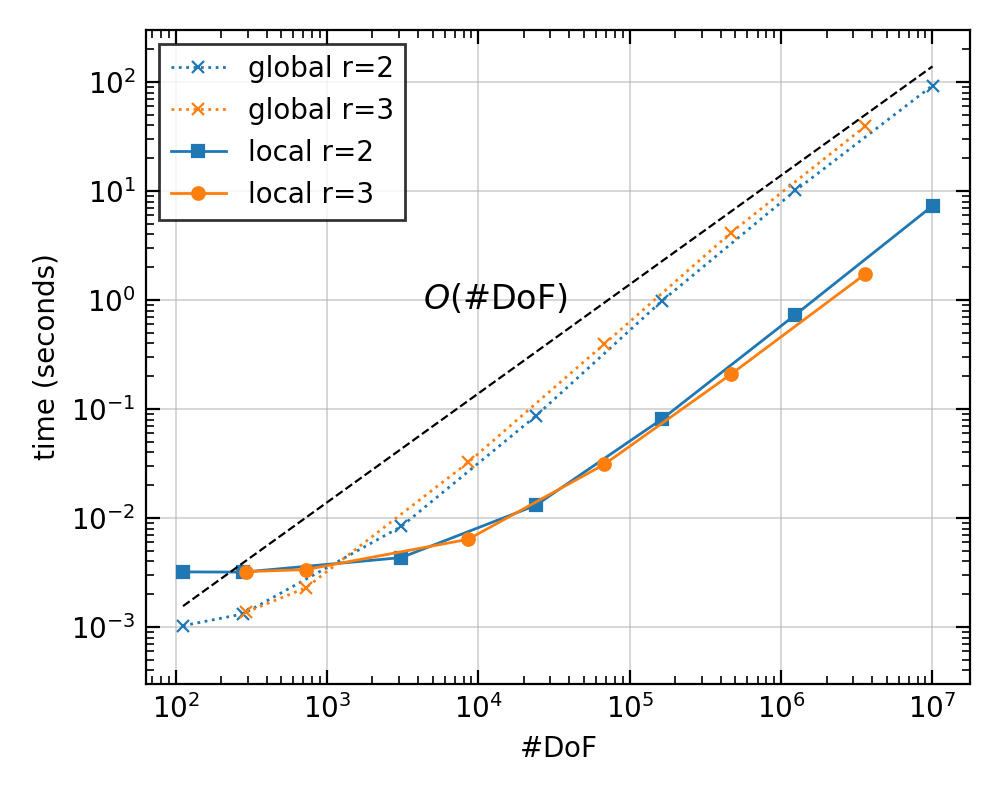}
        \caption{Run time for repeated application.}
        \label{fig:implementation--3d-N1-h}
    \end{subfigure}\hfill
    \begin{subfigure}[t]{0.5\linewidth}
        \centering
        \includegraphics[width=\linewidth]
            {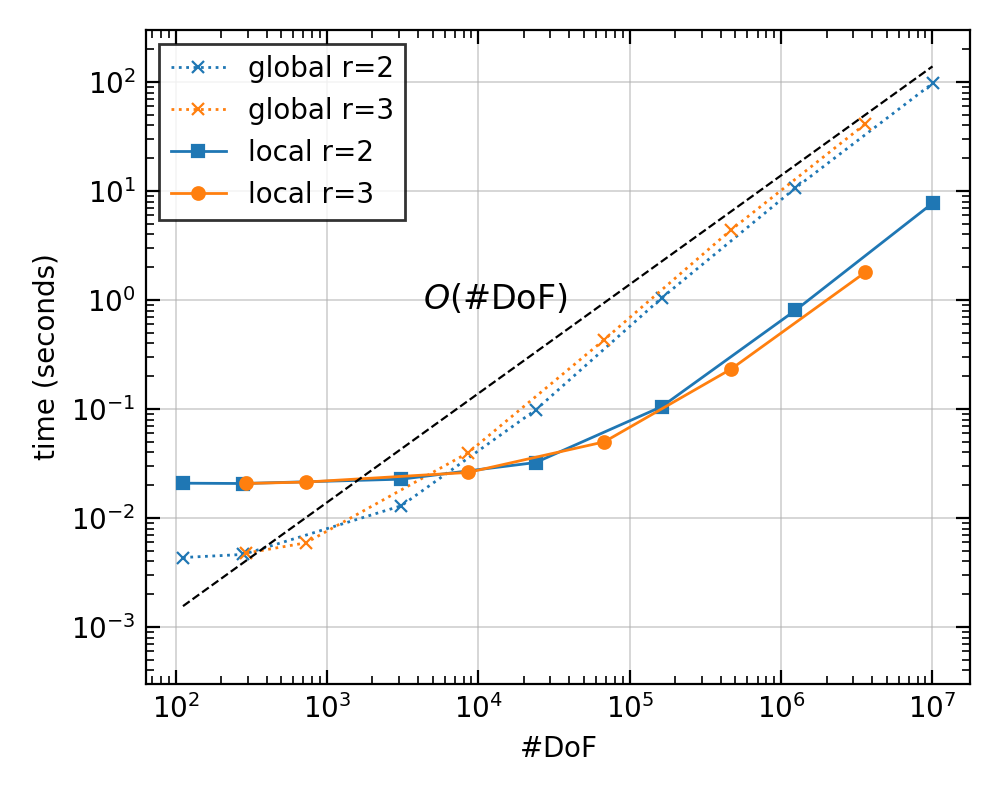}
        \caption{Run time including initialisation.}
        \label{fig:implementation--3d-N1-dof}
    \end{subfigure}
    \caption{Run times for the local discrete Poincar\'e
        operator $\mathscr P_{r,\mathrm{dec}}^1$ acting on the 1st Nédélec space in three dimensions, for polynomial orders
        $r=2,3$. Crosses represent the global baseline. }
    \label{fig:implementation--p1}
\end{figure}

\begin{figure}[h]
    \centering
    \begin{subfigure}[t]{0.5\linewidth}
        \centering
        \includegraphics[width=\linewidth]
            {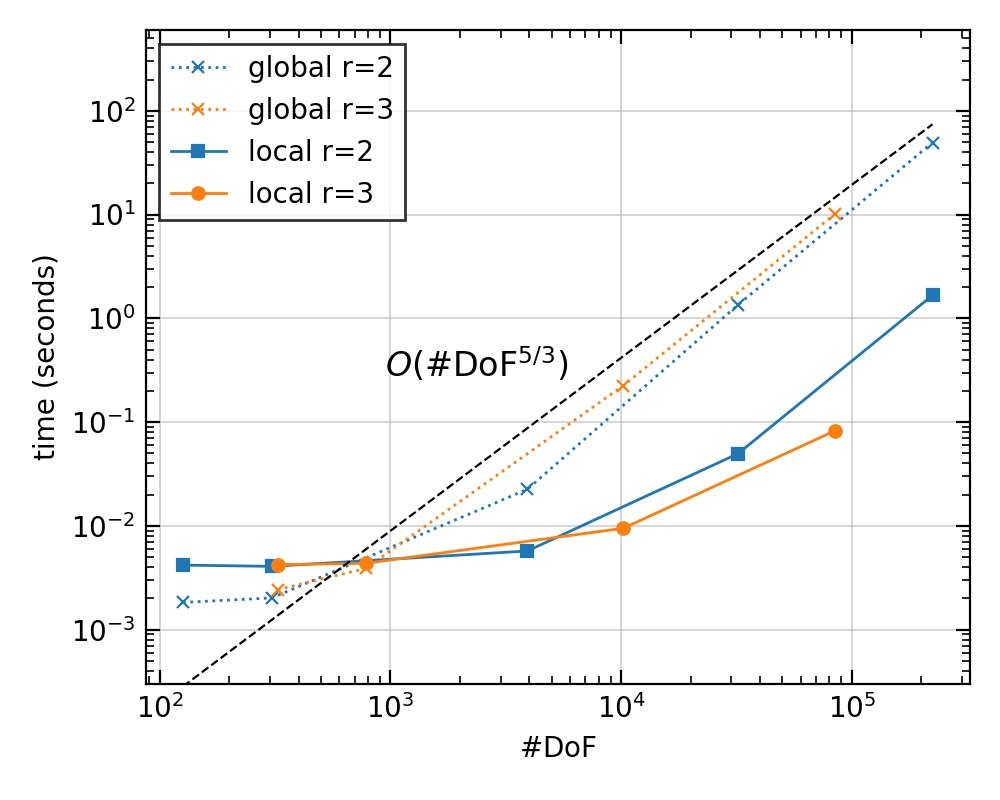}
        \caption{Run time for repeated application.}
        \label{fig:implementation--3d-RT-h}
    \end{subfigure}\hfill
    \begin{subfigure}[t]{0.5\linewidth}
        \centering
        \includegraphics[width=\linewidth]
            {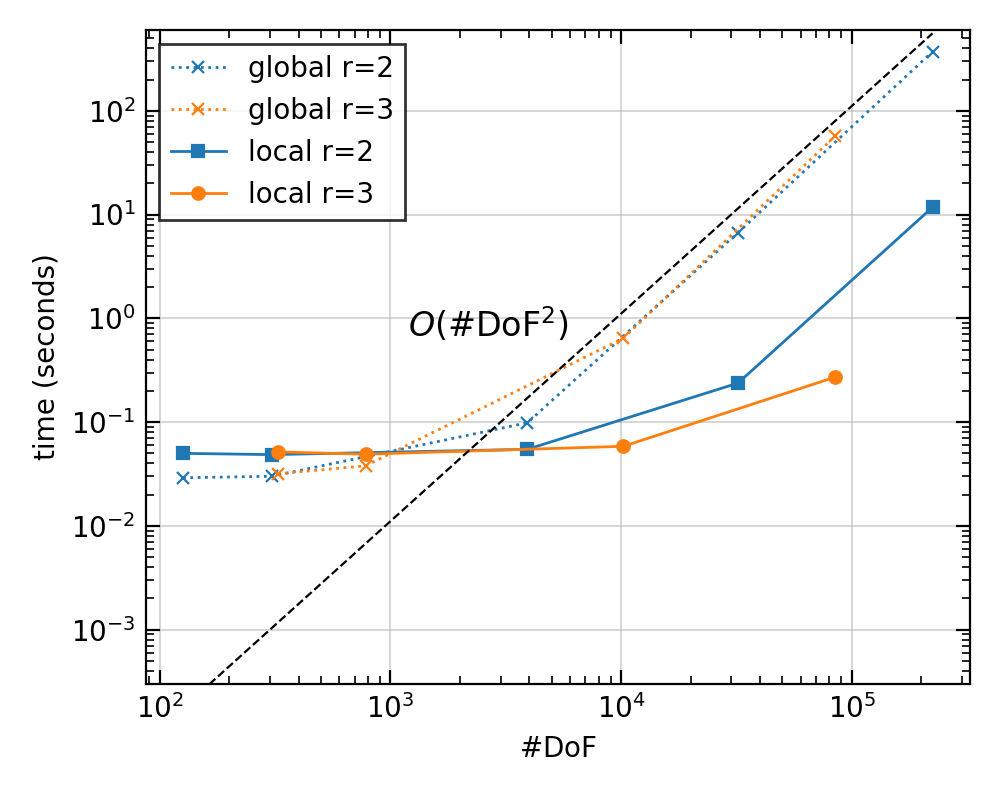}
        \caption{Run time including initialisation.}
        \label{fig:implementation--3d-RT-dof}
    \end{subfigure}
    \caption{Run times for the local discrete Poincar\'e
        operator $\mathscr P_{r,\mathrm{dec}}^2$ acting on the Raviart-Thomas space in three dimensions, for polynomial orders
        $r=2,3$. Crosses represent the global baseline. }
    \label{fig:implementation--p2}
\end{figure}

\begin{figure}[h]
    \centering
    \begin{subfigure}[t]{0.5\linewidth}
        \centering
        \includegraphics[width=\linewidth]
            {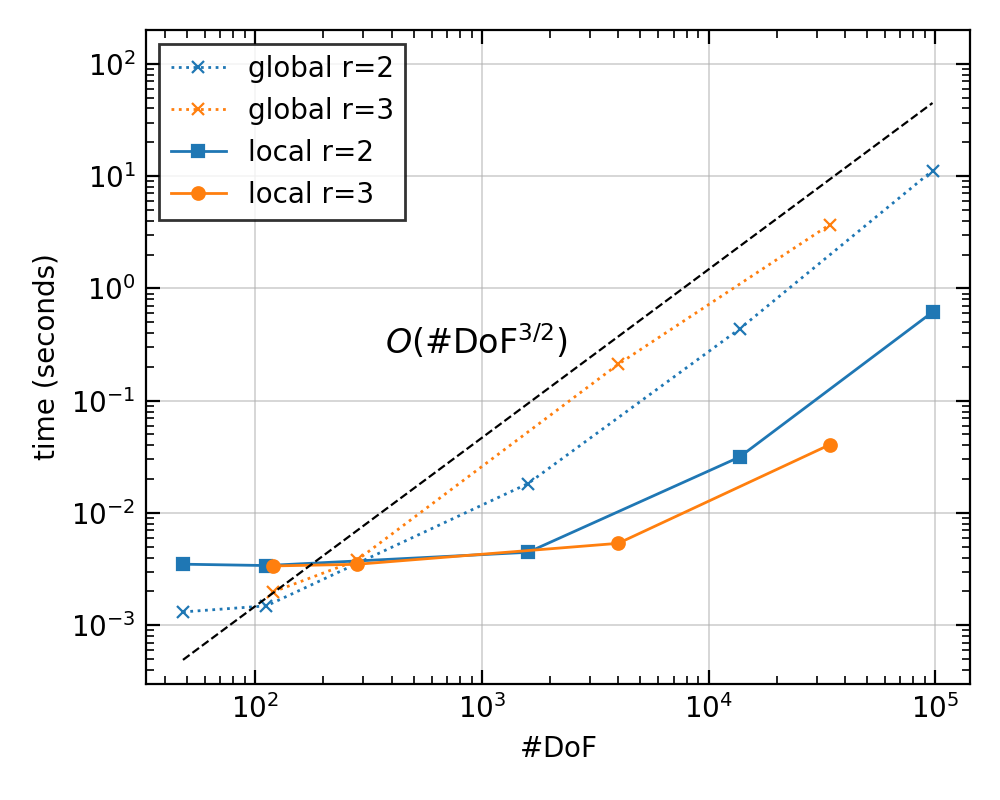}
        \caption{Run time for repeated application.}
        \label{fig:implementation--3d-DG-h}
    \end{subfigure}\hfill
    \begin{subfigure}[t]{0.5\linewidth}
        \centering
        \includegraphics[width=\linewidth]
            {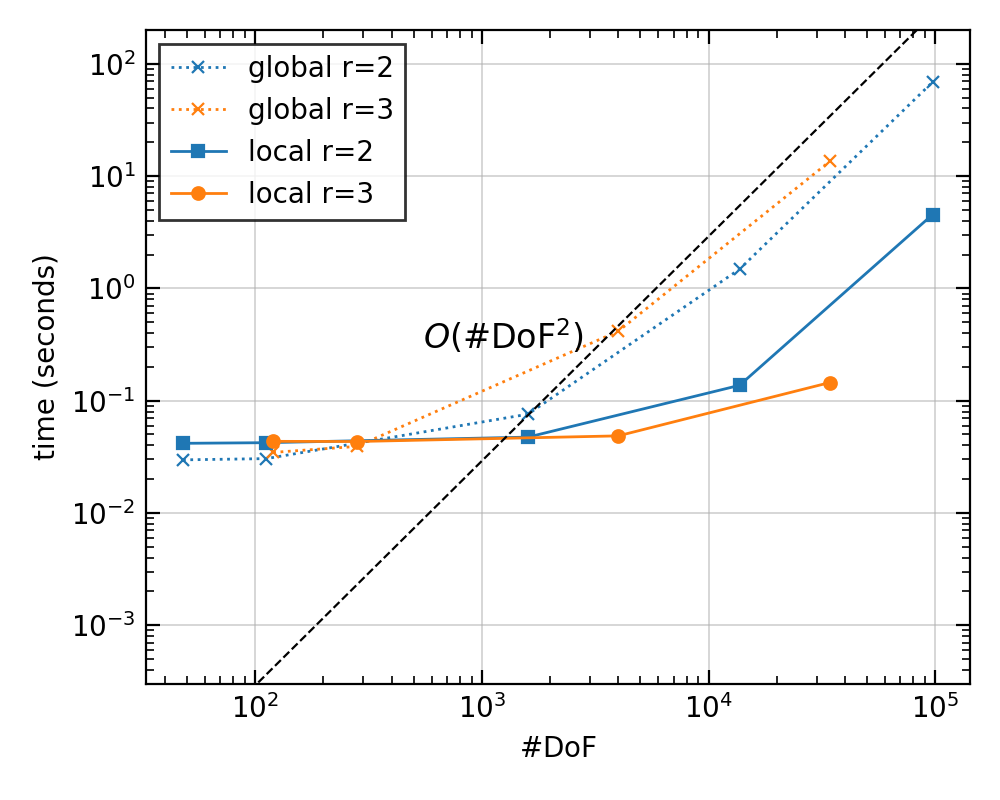}
        \caption{Run time including initialisation.}
        \label{fig:implementation--3d-DG-dof}
    \end{subfigure}
    \caption{Run times for the local discrete Poincar\'e
        operator $\mathscr P_{r,\mathrm{dec}}^3$ acting on the discontinuous Lagrange space in three dimensions, for polynomial orders
        $r=2,3$. Crosses represent the global baseline.}
    \label{fig:implementation--p3}
\end{figure}

The code was ran in Firedrake version 2026.4.2, single-threaded on an Intel Core Ultra 5 335, with 32GB of RAM. \Cref{fig:implementation--p1,fig:implementation--p2,fig:implementation--p3} plot the recorded run times against total degrees of freedom. Each data point is the mean solve time of at least 9 solves, spanning at least 1 second. Solves faster than 100ms are batched, and the solve time is then approximated by the batch average. To test for repeated solves we pre-assemble all the matrices and initialise the relevant Firedrake objects, including the factorised solver. The time needed for this is then excluded from the final plot. For the full-setup test we reassemble and re-initialise everything on each pass. For sufficiently large $\#\mathrm{DoF}$ we see an increased efficiency with the local method in both metrics. 

These are not necessarily the optimal solvers. However, having a better setup for the global problem is equally beneficial for the local problem, since the bulk of the work is solving the lowest-order system. We also want to make clear that the dashed lines are purely for visual reference, and do not necessarily reflect a theoretical time complexity.  

\newpage

\section{Finite element systems}
\label{section:fes}

We conclude by discussing the potential generalizations of the above constructions to other FEEC complexes. This serves two purposes. Firstly, it clarifies the underlying structures that have enabled the above constructions. Secondly, it suggests a path to analogous constructions for more complicated discretizations, like those of the Stokes complex, the elasticity stress complex, and the elasticity strain complex. The appropriate framework is that of finite element systems (FES) \cite{christiansen-rapetti:15,christiansen-hu:22}. It makes for a convenient language, since it exposes the underlying algebraic structure. Here we keep to a surface level account, and refer to \cite{christiansen-rapetti:15} for a proper treatment.

\begin{definition}[Finite element system]
\label{def:fes--system}
A \emph{finite element system} $E^\bullet$ on $\mathcal T_h$ is a family of finite dimensional spaces $E^k(\sigma)\subset H\Lambda^k(\sigma)$, one for each $\sigma\in\mathcal T_h$ and $0\leq k\leq \dim \sigma$, that is stable under restriction and the exterior derivative. In the sense that, for each $\sigma\in\mathcal T_h$, it holds that $\tr_\tau E^k(\sigma)\subset E^k(\tau)$ for all $\tau\triangleleft\sigma$, and $dE^k(\sigma)\subset E^{k+1}(\sigma)$. 
\end{definition}

\Cref{prop:prelim--extder-polynomial,prop:prelim--trace-polynomial} say that $\mathcal P_r\Lambda^\bullet$ and $\mathcal P_r^-\Lambda^\bullet$ are finite element systems. A less standard example is the remainder family $\mathcal M_r^-\Lambda^\bullet$ of \eqref{eq:lc--remainder-spaces}, which is a finite element system by \cref{lem:lc--interpolant}. The structure is set up so that we can define global spaces by gluing. For a sub-mesh $\mathcal S_h\subset \mathcal T_h$, we define its global space by
\begin{equation}
\label{eq:fes--global}
E^k(\mathcal S_h)\coloneqq \set[\bigg]{u\in\bigoplus_{\sigma\in\mathcal S_h}E^k(\sigma):
\tr_{\sigma\cap\sigma'}u_\sigma=\tr_{\sigma\cap\sigma'}u_{\sigma'}\quad\forall\sigma,\sigma'\in\mathcal S_h}.
\end{equation}
In the context of FES, one is particularly interested in $E^k(\mathcal T_h)$ and $E^k(\partial \sigma)$. Note that the construction of $E^k(\mathcal T_h)$ is the same as in \cref{eq:prelim--global-assembly}, except that we are adding in some redundant information on the lower-dimensional cells.

\subsection{Degrees of freedom}

Going forward we fix a finite element system $E^\bullet$. Furthermore, we  suppose we have identified some \emph{dual} finite element system $Q^\bullet$, such that on every sub-mesh $\mathcal S_h\subset \mathcal T_h$ the set of dof given by
\begin{equation}
\label{eq:fes--dof}
u\mapsto\int_\tau\tr_\tau u\wedge q,
\qquad q\in Q^{\ell-k}(\tau),
\end{equation}
for all $\tau\in\mathcal S_h^\ell$ and $\ell \geq k$, is unisolvent on $E^k(\mathcal S_h)$. We also make the assumption that $E^0(\sigma)$ and $Q^0(\sigma)$ both contain the constant functions. This generalizes the notion of canonical dof (recall \cref{prop:prelim--dof-full-polynomial,prop:prelim--dof-trimmed-polynomial}). Indeed, the test forms that appear in \cref{prop:prelim--dof-trimmed-polynomial} as canonical dof for $\mathcal P_r^-\Lambda^k$ are given by
\begin{equation*}
Q^{\dim\sigma-k}(\sigma)=\mathcal P_{r-\dim\sigma +k-1}\Lambda^{\dim \sigma-k}(\sigma).
\end{equation*}
In \cref{section:prelim--global-spaces} we made the observation that this family, as well as the analogous family for $\mathcal P_r\Lambda^\bullet$, is itself stable under the trace and the exterior derivative. This is equivalent to saying that $Q^\bullet$ is a finite element system. These kinds of structures were originally studied in \cite{christiansen:11}, and then later in \cite{christiansen-rapetti:15}. 

We now want to know whether the local construction of \cref{section:local-construction} translates to the general setting. This requires the notion of a lowest order space, by which we mean some FES $\mathfrak W^\bullet_E$ subordinate to $E^\bullet$, which is analogous to the Whitney forms. Such spaces were obtained in \cite{christiansen-rapetti:15} through the decomposition by harmonic degrees of freedom. Here we consider the analogous construction with respect to dual systems. Given a base point $p_\sigma\in \sigma$, and assuming $Q^0(\sigma)$ is sufficiently regular, we have the direct sum decomposition
\begin{equation}
Q^0(\sigma) = \bigl\{q\in Q^0(\sigma):q(p_\sigma)=0\bigr\} \,\oplus\, \mathbb R,
\end{equation}
By zeroing all the dof except those corresponding to $\mathbb R$ we obtain $\mathfrak W^\bullet_E$. The local space $\mathfrak W^k_E(\sigma)$ consists of all forms $u\in E^k(\sigma)$, such that for all $\tau\triangleleft \sigma$ and $q\in Q^{\dim \tau - k}(\tau)$, it holds
\begin{equation}
\int_\tau u \wedge q = 
\begin{cases}
    0 & \text{if $\dim \tau > k$}, \\
    q(p_\tau)\int_\tau u & \text{if $\dim \tau=k$}.
\end{cases}
\end{equation}
Then $\mathfrak W^k_E$ carries the same set of canonical dof as $\mathcal P_1^-\Lambda^\bullet$
\begin{equation}
u \mapsto \int_\tau \tr_\tau u \qquad \forall \tau\in \mathcal T_h^k,
\end{equation}
and by equating these we obtain a natural isomorphism between $\mathcal P_1^-\Lambda^\bullet$ and $\mathfrak W^k_E$.\footnote{In the sense that each local isomorphism $\Phi_{k,\sigma}:\mathcal P_1^-\Lambda^k(\sigma)\to \mathfrak W^k_E(\sigma)$, and its inverse, commute with restriction and the exterior derivative.} In particular, this lets us define a canonical interpolator $\tilde I^k_\sigma:E^k(\sigma)\to \mathfrak W^k_E(\sigma)$, in the same way we did in \cref{section:lc--decomposition}. Thus we recover the decomposition $E^k(\sigma)= \mathfrak W^k_E(\sigma)\oplus \mathcal M_E^k(\sigma)$, with $\mathcal M_E^k(\sigma)\coloneqq \ker \tilde I^k_\sigma$.

\begin{lemma}
\label{lemma:fes--cohomology}
If $E^\bullet$ can be equipped with a dual system, then the cohomology of $E^\bullet(\mathcal T_h)$ is isomorphic to the de Rham cohomology of $\Omega$.
\end{lemma}

\begin{proof}
Since $Q^\bullet$ is a FES, it follows from Stokes' theorem
\begin{equation}
    \int_\sigma du\wedge q = \sum_{\tau\triangleleft \sigma}\int_{\tau} \tr_\tau u\wedge \tr_\tau q + (-1)^{k+1}\int_\sigma u\wedge dq
\end{equation}
that the canonical interpolant onto $E^\bullet$, as induced by $Q^\bullet$, commutes with the exterior derivative (see \cite[Equation 24]{christiansen-rapetti:15}). Then $E^\bullet$ is ``compatible'' in the sense of \cite[Proposition 2.8]{christiansen-rapetti:15}, and the statement follows from standard results (see e.g. \cite[Proposition 5.16]{christiansen:11}).
\end{proof}

\begin{proposition}
\label{prop:fes--remainder-cohomology}
Let $E^\bullet$ be a finite element system equipped a dual system $Q^\bullet$. For every $\sigma\in\mathcal T_h^m$, the complexes of local remainders are exact
\begin{subequations}
\begin{equation}
\label{eq:fes--local-exact}
0\longrightarrow{\mathcal M}_E^0(\sigma)
\xrightarrow d\cdots\xrightarrow d
{\mathcal M}_E^m(\sigma)
\longrightarrow0;
\end{equation}
\begin{equation}
\label{eq:fes--local-exact-bubble}
0\longrightarrow\mathring{\mathcal M}_E^0(\sigma)
\xrightarrow d\cdots\xrightarrow d
\mathring{\mathcal M}_E^m(\sigma)
\longrightarrow0.
\end{equation}
\end{subequations}
\end{proposition}

\begin{proof}

By \cref{lemma:fes--cohomology} the cohomology of the local complex $E^\bullet(\sigma)$ is isomorphic to that of $\mathcal P_r^-\Lambda^\bullet(\sigma)$. By \cite[Proposition 2.3]{christiansen-rapetti:15} the same holds for the cohomology of the local bubbles $\mathring{E}^\bullet(\sigma)$ with respect to $\mathring{\mathcal P}_r^-\Lambda^\bullet(\sigma)$. The rest of the proof is the same as in \cref{section:lc--decomposition}.
\end{proof}

From here we are free to go about the construction of \cref{section:local-construction} in much the same way. The construction of \cref{section:lc--Bk-construction} only uses the local spaces on each $\sigma$ and its faces, together with their traces, degrees of freedom and inner products, so it applies verbatim; only the reference simplex variant of \cref{rem:lc--cost} additionally requires the local spaces to be defined in terms of the reference simplex. 

We also check that the construction of \cref{section:cone-extension} carries over. 

\begin{proposition}
\label{prop:fes--poincare}
Let $\Omega$ be Lipschitz contractible in the sense of \cref{def:conext--Lipschitz-contraction}, and let $E^\bullet$ be a finite element system equipped with a dual system $Q^\bullet$. If each element of $E^\bullet(\mathcal T_h)$ is bounded near $p_0$, then
\begin{equation}
\label{eq:fes--poincare-def}
\int_\sigma \mathscr P^{k+1}u\wedge q \coloneqq \int_{\cone(\sigma)} \Phi_\sigma^*u\wedge \pr^*_\sigma q
\qquad \forall q\in Q^{\dim\sigma-k}(\sigma),\ \sigma\in\mathcal T_h ,
\end{equation}
defines a Poincar\'e operator $\mathscr P^\bullet$ on $E^\bullet(\mathcal T_h)$ in the sense of \cref{def:prelim--poincare}.
\end{proposition}

\begin{proof}
Unisolvence of \eqref{eq:fes--dof} makes $\mathscr P^{k+1}u$ a well defined element of $E^k(\mathcal T_h)$. The proof of the homotopy identity then follows as noted in \cref{section:conext--general-construction}. The polynomial structure only entered in two places: through stability of the test forms under $\tr_\tau$, used to pass from $\sigma$ to its facets in \eqref{eq:conext--homotopy-proof-term11}, and through their stability under $d$ to handle $dq$. These are the two requirements in \cref{def:fes--system} applied to $Q^\bullet$.
\end{proof}

\subsection{Comment on discretizations of Stokes and elasticity}

A large body of work has been developed to construct finite element subcomplexes of the Stokes and elasticity complexes; see, e.g. \cite{JLMN:17, Guzman-Neilan:18, Chen-Huang:22, Chen-Huang:24, Chen-Huang:26, CGGH:24}. In particular, the works \cite{CGGH:24, Chen-Huang:26} use the Bernstein-Gelfand-Gelfand (BGG) construction to develop finite elements for elasticity. Bounded continuous-level Poincar\'e operators have been developed for the BGG complexes in \cite{Cap-Hu:23}. A natural question is thus whether analagous discrete Poincar\'e operators can be developed for these complexes.

The framework of finite element systems has been generalized to also cover the Stokes and elasticity complexes \cite{christiansen-hu:22,christiansen-hu:18}. The treatment is analogous to what we've done so far for the de Rham complex, except that the spaces on lower-dimensional cells are bulked to carry more information, for example in order to accommodate higher regularities. This in turn makes the restriction and differential maps more complicated, especially as you go to discretisations in higher dimension. Cell-wise decompositions into local complexes with commuting restriction maps have been identified for the two dimensional Clough-Tocher and Johnson-Mercier finite elements \cite{christiansen-hu:22}. Moreover, such decompositions carry lowest order spaces, analogous to Whitney forms \cite[Remark 18]{christiansen-hu:22}. Since the algebra of generalized FES remains largely unchanged, it seems to us reasonable to expect the construction of \cref{section:local-construction} to carry over.

\section*{Declaration of Generative AI use}

During the preparation of this work the authors used Claude Opus and ChatGPT to assist with coding, literature review, figure generation, and proofreading. The authors have reviewed and edited the content as needed and take full responsibility for the content of the published article.

\printbibliography{}

\end{document}